\documentclass[11pt]{article}
\usepackage{amsmath,amssymb,latexsym,float,epsfig,subfigure}
\usepackage{latexsym}
\usepackage{natbib}
\usepackage{amsthm}
\usepackage{cases}
\usepackage{epsfig}
\usepackage{subfigure}
\usepackage{scalefnt}
\usepackage{booktabs}
\usepackage[english]{babel}
\usepackage{graphicx}
\usepackage{titlesec}
\usepackage{epsfig}
\usepackage{enumerate}
\usepackage{caption}
\usepackage{subfigure}
\usepackage{booktabs}
\usepackage[english]{babel}
\usepackage{graphicx}
\usepackage{multirow}
\usepackage{rotating}
\usepackage{fancyhdr}
\usepackage{setspace}
\usepackage[section]{placeins}
\usepackage{xcolor}

\usepackage{breqn}

\usepackage{array, makecell}
\usepackage{pseudocode}
\usepackage{algorithmicx}
\usepackage{algorithm}
\usepackage{url} 
\usepackage{appendix}
\usepackage{ragged2e}
\usepackage{physics}

\usepackage{chngpage}

\usepackage[nocomma]{optidef}
\usepackage{mathtools,zref-savepos}
\newtagform{roman}[]()

\topmargin by -\headsep \textheight 9in \oddsidemargin 0pt
\evensidemargin \oddsidemargin \marginparwidth 0.5in \textwidth
6.5in
\begin{document}
\title{Word-of-Mouth on Action: Analysis of Operational Decisions When Customers are Resentful}

%
\author{Bahar \c{C}avdar \thanks{bcavdar@tamu.edu} \and Nesim Erkip \thanks{nesim@bilkent.edu.tr}}
\date{
    $^*$\small{Department of Engineering Technology and Industrial Distribution, Texas A\&M University, College Station, TX, USA }\\%
    $^\dagger$Department of Industrial Engineering, Bilkent University, Ankara, Turkey\\[2ex]%
}

\maketitle

\begin{abstract}
Word-of-Mouth (WoM) communication, via online reviews, plays an important role in customers' purchasing decisions. As such, retailers must consider the impact of WoM to manage customer perceptions and future demand. In this paper, we consider an online shopping system with premium and regular customers. Building on the behavioral and operations management literature, we model customer preferences based on the perceived service quality indicated by WoM and integrate this into the retailer's operational problem to determine a shipment policy. First, we study the e-tailer's problem when she has no knowledge of WoM, and only reacts to the changes in demand. We analyze the long-term behavior of customer demand and show that potential market size and customer sensitivity are the key parameters determining this behavior. Then, we build a model to integrate knowledge of WoM into operational decision making and partially characterize the optimal solution. We show that (i) beating the competition in the market may not always benefit the company under strict operational constraints as it can result in undesired customer switching behavior, (ii) relaxations in operational constraints may hurt profitability due to the associated difficulties of managing perceptions, and (iii) seeking a stationary policy can lead to suboptimal solutions, therefore cyclic policies should also be considered when appropriate. 
\end{abstract}

\section{Introduction}

Advances in online retailing have introduced new modes of sales, an ease of collecting customer-related data, and new ways for sellers to interact with customers both pre-purchase and post-purchase. Additionally, customers have new channels to communicate with each other as a form of electronic Word-of-Mouth (WoM). One of them is online reviews; customers share their experiences, ask questions, and make complaints online. Studies show  that 93\% of customers trust online reviews for their purchase decisions (\cite{qualtrics}).

WoM can provide information on a variety of service attributes including service quality. By changing the quality, the company can affect the impression of the service, and consequently consumer choices. This creates a feedback mechanism between operational decisions and demand. Our goal is to investigate environments where operational decisions impact perceived service quality, and perceived service quality impacts demand in return. We analyze an environment where a company offers different delivery durations (as different services) to supply customers, and demand is endogenous to operational decisions through WoM.

Offering different service options to increase customer base and profits is a common practice in different industries. Hotels (e.g., Hilton) run loyalty programs to attract returning customers by offering free-of-charge upgrades, and discounted additional services. Members-only vacation clubs (e.g., Happimag) give higher priority to members over regular customers in booking and other services, and offer them better prices. Retailers offer premium services for a charge to cover shipping expenses and/or to expedite a delivery, as well as some products or services that may not be available elsewhere at the time (e.g., Amazon). As seen from these examples, it is common to create a deliberate difference in service quality to attract customers who are willing to pay more.

In the aforementioned service environments, customers can be classified in two groups, namely premium and regular customers. Premium customers pay membership fees for a guaranteed high-quality and faster service. Regular customers, on the other hand, wait at least as much as premium customers, but likely more. When customers are deciding on whether to choose premium or regular service, they analyze prior reviews of the services, i.e., WoM. If the perceived difference in service quality between premium and regular services is high, there is a greater incentive to purchase the premium service. However, if the perceived difference is low, then there is little benefit to purchase and fewer customers will do so. 

We study the impact of WoM communication in the setting of an online shopping delivery system where an ``e-tailer" offers a homogeneous product to regular and premium customers. The e-tailer's operational strategy dictates WoM reviews thereby influencing future customer demand. In turn, the e-tailer adjusts the operations policy to respond the new demand creating a feedback between operational decisions and customer service selection.

 \subsection{Motivation and Background}
 
Our work is mainly motivated by the studies addressing endogenous demand structures in the Operations Management (OM) literature. An early work that motivates this study is by \cite{bitran2008om}. They form a connection between operational policies and customer behavior, and emphasize a need to bridge the gap between mathematical representation of the impact of service delivery duration on customer satisfaction and the current knowledge of human behavior from psychological and marketing literature. \cite{netessine2009consumer} also emphasize the importance of endogenous demand especially for operational decision making, and state that it is necessary to develop such models. Addressing such an issue, we operationalize these conceptual ideas, and model the relation between demand and operational decisions. \cite{liu2019channel} consider endogenous demand where customers are sensitive to the lead time and and analyzes channel strategies (direct, traditional or hybrid) and online pay-on-delivery services.
We focus on determining the right operational policy to manage customer perceptions for an e-tailer. We also  consider the idea of market expansion by including premium service. Specifically, we follow the work by \cite{kumar2010exploiting} for the use of an equilibrium analysis to compute the size of the additional market that can be captured, and \cite{robinson2016appropriate} for the use of an optimization approach to determine the profitable size of the additional market.

\subsection{Contributions}
Building on the literature at the intersection of behavioral and operational issues, we study an e-tailer's shipment policy problem when customer demand is affected by operational decisions through WoM. Our main contributions are based on leveraging knowledge of WoM in operational decision making, and providing tools for a better management of customer perceptions in a profit-oriented company. More specifically, we build a model for customer demand as a WoM-initiated feedback function of operations, develop a stable strategy for operations, analyze the sensitivity of profitability and operational decisions to key parameters under WoM, and analyze the value of knowledge of WoM.

Our results provide insights for industries where customer switching between different service providers or services is common. Main insights can be summarized as follows:
\begin{itemize}
\item[] Beating the competition in the market and providing a better service than other companies may not always benefit the company under operational constraints. Indeed, it can create a sudden increase in reputation leading to an undesired demand inflow that is beyond the service capacity. As a result, the good reputation cannot be sustained, and the profitability over time will be subject to market dynamics rather than being under control.

\item[] Investing in relaxing operational constraints expecting that it will increase profitability may not turn out to be so; in fact, the relaxations can hurt profitability due to the increased difficulty of managing customer perceptions when these constraints are associated with customer perceptions.

\item[] Although cyclic operational policies are not desirable in general, we show cases where they can be more profitable than stationary policies. Hence, non-stationary equilibrium policies should be considered as well.
\end{itemize}

These insights also illuminate the examples of company bankruptcies due to failing to understand endogenous demand structures in early 2000's (e.g., bankruptcy of Charming Charlie) as it will be discussed later.

The rest of the paper is organized as follows. In Section \ref{Sec:Lit}, we review the relevant literature. In Section \ref{Sec:Analysis}, we present our models and theoretical results characterizing the optimal solutions. We present and discuss some numerical results in Section \ref{Sec:Comp}, and conclude with Section \ref{Sec:Conc}.

\section{Literature Review}
\label{Sec:Lit}

To understand the relation between operational decisions and demand, we draw on interdisciplinary literature on WoM and customer loyalty. We also review the operational aspects of e-tailers. Finally, we review operations management literature on marketing and operations interface. We present an overview of these topics in relation to our work, and position our paper with respect to the literature.

\subsection{Demand Structure: WoM and Loyalty Concerns}

We analyze the related literature to develop two main perspectives on the following issues: (i) how demand is affected by WoM, and (ii) the relation between the perceived service quality, and loyalty and switching issues. For the latter, we especially focus on the effect of waiting time information which is a main determinant for the service quality perception for online retailing (\cite{rao2011failure}).

In WoM, customer statements are considered to be a way of communication. \cite{hennig2004electronic} define customer statements as ``any positive or negative statements made by potential, actual or former customers about a product or company, which is made available to a multitude of people and institutions via internet''. There are many different WoM (more specifically electronic WoM) communication channels. Some examples are social network structures (\cite{luis2013drivers}), online reviews (\cite{park2008effects}) and micro blogging (\cite{jansen2009twitter}), etc. \cite{ramanathan2017role} study social media reviews as a reflection of customers' perception of service. \cite{besbes2018information} claim that use of online ratings to learn the quality of a product in an unbiased way is possible under some general assumptions. 

\cite{bayus1985word} is one of the first to mention a connection between WoM and marketing operations. He considers a feedback mechanism between WoM and marketing efforts, and proposes a model that assumes WoM and marketing efforts are both affecting and affected by sales. 

There is a strong evidence showing that positive evaluations increase loyalty/demand for a service. Therefore, we also review the literature concerning loyalty and customer switching behavior. \cite{heim2001operational} investigate the drivers of customer loyalty in electronic retailing. One of their empirical findings indicate that operational performance is a significant driver in customer loyalty. \cite{rao2011electronic} present empirical studies that relate the service quality with future demand. \cite{heim2001operational} and \cite{kumar2011impact} show that customer loyalty partly depends on operational performance. 

Switching is generally associated with the perception of the difference between two services, e.g., if the difference between premium and regular services is large, customers are more likely to switch to the premium service. To model the switching behavior, we use the concept of regret from the seminal work of \cite{bell1982regret}, which is based on defining a disutility of a possible regret related to the outcome. The idea of using disutility can also be followed in the works of \cite{roels2013optimal}, \cite{veeraraghavan2011herding}, and \cite{rong2014role}. They show that sunk costs and service quality of an initial choice and new choice affect the level of regret felt regarding the consumer's initial choice through empirical studies. 

The idea of endogenous demand structure displayed by the WoM literature is also supported by the OM literature. There are many papers considering a relation between demand and service quality and price. \cite{veeraraghavan2011herding}, \cite{liu2013pricing} and \cite{afeche2017customer} are some examples. We use this idea while modeling the impact of relative difference between different types of services on customer demand. \cite{huang2021influence} is one study that is closely related to ours. They model the impact of online reviews on customers' views about the actual service quality. They show that when the capacity is flexible, companies may prefer not to adopt online reviews since they can easily adjust the fluctuations in demand. However, companies with limited capacity should adopt online reviews to allow customers to have accurate estimates about the quality of their service to better manage the demand.

Based on the above literature review, we assume that customers can observe WoM, and form reactions through a response function. Additionally, we consider that there is a relation between service performance and loyalty, and switching behavior. Bringing these together, we model the impact of relative difference between different types of services on premium customer demand.

\subsection{Conducting Operations for E-tailers }

The importance of operational performance for e-tailers is emphasized by several studies (e.g., \cite{kumar2011impact} and \cite{fisher2016value}). When customers are informed about the waiting time (usually maximum delivery time is declared), they feel less irritated because they expect it (\cite{taylor1994effects}), and as a result they consider it to be more acceptable and give better service evaluation (\cite{hui1996tell}). Similarly, providing waiting time information decreases the perception of the actual waiting (\cite{ahmadi1984effects}), and gives customers a feeling of higher control (\cite{hui1996does}). \cite{nie2000waiting} studies social and psychological factors affecting the perception of waiting on the customer end and discusses strategies service providers can utilize to reduce the perceived waiting time.

Following these papers, we assume that when provided with the maximum waiting time information, regular customers accept it as the expected delivery duration. It is important for them that this promise is kept, and their satisfaction increases with the speed of the delivery. On the other hand, premium customers mostly care about the differential speed between the two types of services. In our model, we relate the service level of regular customers to the demand generated by premium customers using regret concept: the smaller the gap between premium and regular service quality, the lower the demand for the premium service. 

Operations strategy for an e-tailer has multiple components. One component of an e-tailing environment is related to how operations are carried out. Vehicle consolidation is one of the components (for related recent work on shipment consolidation see \cite{mutlu2010integrated}, \cite{ccapar2013joint}, \cite{kaya2013coordinated}, \cite{ulku2012optimal}, \cite{ulku2012modelling}). The other component is order fulfillment approach (see for instance \cite{netessine2006supply} for drop shipping idea, \cite{chiang2005managing} for a multi-echelon structure of fulfillment centers). A final component of an operational strategy is regarding strategic lost-sales (\cite{anand2011quality}). 


We also review the studies on the cost/profit structure relevant to the e-tailing environment we consider. The premium membership fee (as a fixed cost) has a number of important roles in customer purchase behavior. We refer to the work of \cite{arkes1985psychology}, \cite{dick1998impact} and \cite{wang2010sunk} to structure the fee related functions, to be described later.


\subsection{Relevant Operations Management Literature}

A stream of OM literature focuses on the interaction between different aspects of the problem. For example, \cite{eliashberg1993marketing} review marketing-production joint decision making. \cite{chan2004coordination} focus on coordination of pricing and inventory decisions. \cite{yano2005coordinated} review coordinated pricing with production and procurement decisions. However, these papers mainly focus on the integration of different operational decisions. There are only a few studies considering WoM as a part of this interaction as we do in this paper. 

When considering a dynamic relation between different aspects of a problem, using an equilibrium analysis is common. \cite{erickson2011differential} presents an equilibrium analysis for marketing-operations interface at a strategic level. \cite{grosset2011goodwill} consider advertising as a management function and analyze the resulting equilibrium demand. \cite{yan2009optimal} consider a new product development case where WoM feedback eventually results in an equilibrium for the potential sales of the product. \cite{tereyaugouglu2012selling} find an equilibrium pricing and production policy relating customer beliefs to service level. 

The papers that are closely related to ours are by \cite{aflaki2013managing} and \cite{afeche2017customer}. The former studies customer retention using loss aversion ideas. They consider the probability of a customer to bring a new customer, which is interpreted as WoM. The latter analyzes the relation between the capacity allocation and customer retention. They consider a negative impact of rejecting a customer on attracting new customers. These papers are significant in terms of modeling the impact of WoM; however, they treat it as a very general information. We model WoM based on different attributes of the service. Our paper considers the quality of regular service as a source of disutility for premium customers, and hence is different from the current literature. We take a similar approach to the literature to leverage knowledge of WoM and seek stationary equilibrium policies. Finally, we carry out our analysis in a deterministic environment to simplify the exposition.

\section{Problem Definition and the Analysis  }

\label{Sec:Analysis}

In the first part of this section, we present the introductory elements of our problem definition. We then analyze the e-tailer's policy when she does not have any knowledge of WoM. We conclude the section with a model to jointly optimize operations and manage WoM-based customer behavior to maximize her average profit based on the assumption that the e-tailer knows how WoM impacts demand. 

\subsection{Description of Environment and Modeling Framework}
\label{Sec:NoKnow}
In this section, we introduce the basics of the operational environment and the elements of our models, the impact of premium membership fee, and the structure to represent WoM and its influence on premium customer demand.

\subsubsection{E-tailer Operations and Policy Structure}
\label{Sec:operations}
We consider an e-tailer offering a single product to two types of customers concentrated in one geographical area: premium customers who are charged a membership fee to enjoy a guaranteed service with negligible delivery time and regular customers. The service term in our study is defined as the time to deliver an order.

When a regular customer gives an order, the company can wait up to $\tau$ days to fulfill the order, where $\tau$  can be defined as the maximum waiting time announced for the delivery of regular customer orders, i.e., maximum delivery duration. 
Following the loyalty notions discussed in the literature, we consider regular customer demand rate to be constant assuming that non-member orders are drawn from a much bigger population and are not sensitive to the membership terms.

We consider a two-echelon inventory system that involves a fulfillment center and a small local depot to serve different customer types from each echelon. The fulfillment center is the main shipment source of regular customer orders, and supplies the depot. The depot is located geographically close to the location of the anticipated demand from premium customers. Maintaining a central location as an inventory point incurs additional holding cost of $h$ per unit per time, yet it enables faster delivery. The depot predominantly serves premium customers but can also satisfy regular demand. The option of satisfying some regular customer orders immediately from the depot as opposed to sending them with the next shipment batch from the fulfillment center is not considered in most of the previous studies. This setting is a version of drop-shipping idea described by \cite{netessine2006supply} and the two-echelon inventory system mentioned by \cite{chiang2005managing}.

The above structure allows us to consolidate the shipments for fulfilling the orders of regular customers with the replenishment of the local depot inventory. 
This is in line with Amazon's Anticipatory Shipping, where items are sent to the destination geographical area before purchase. While it creates a form of inventory, it is still a small amount (compared to classical inventory replenishment) that can be carried by trucks used for order delivery. We assume unlimited truck capacity with a fixed shipment cost $K$.  As in many real-life applications, we do not charge premium customers any extra shipping cost. To increase flexibility and profit, we strategically allow lost sales possibility for regular customers only (\cite{cheung2016submodular}), whereas we satisfy all premium customer demand as promised. We assume that we receive the same revenue $r$ per sale from any customer type. Finally, we assume that demand rate is deterministic for both customer types to simplify the exposition. For premium customer demand, we allow WoM to be one of the drivers as it will be described in Section \ref{Sec:Know}.

We do not impose a specific inventory replenishment rule. However, following the literature we seek a stationary shipment policy with a shipment cycle of time length $T$. Note that in each shipment cycle, i.e., replenishment cycle, a truck leaves the fulfillment center to deliver regular customers' orders and replenish the local depot inventory.

We can further characterize a shipment cycle. First recall that all premium demand during a shipment cycle must be satisfied. As the cycles are repeating, in each cycle there may be different phases with each phase representing a different treatment for regular customer demand. 
\vspace{-0.1in}
\begin{itemize}
\item[] Phase 1 (fast service): Regular customer orders are sent from the depot immediately after orders are received. 
\item[] Phase 2 (no service): Regular customer demand is strategically lost. \item[] Phase 3 (regular service): Regular customers receive the usual service, i.e., their orders are sent from the fulfillment center within $\tau$ days.
\end{itemize}

Note that in all phases, premium orders are sent immediately from the depot. Finding a \emph{shipment policy} for the e-tailer is to determine the length of the shipment cycle $T$ and decide on how each phase is placed in the shipment cycle. The shipment policy also determines the required inventory level for the local depot. The e-tailer's problem is to find an optimal shipment policy that maximizes the average profit per unit time.

\begin{figure}[H]
	\begin{center}
		\includegraphics[ width=0.7\textwidth]{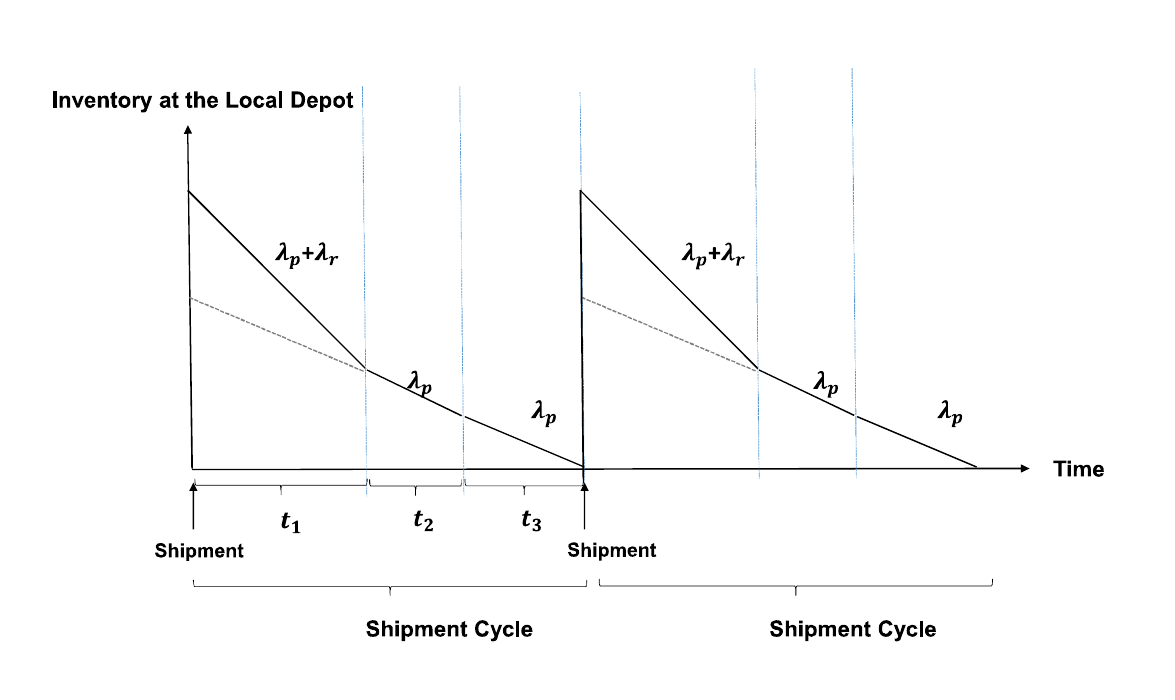}
		\caption{Phases of a shipment cycle.}
		\label{Fig:cycle}
	\end{center}
\end{figure}

If we consider an e-tailer solving this problem based on current premium and regular customer demand without any knowledge of WoM, we can further characterize the optimal shipment policy: In Phase 1, we carry additional inventory in the local depot; hence, Phase 1 must be the first phase of a shipment cycle to avoid additional inventory cost. The duration between the start of Phase 3 and the next shipment is constrained by the maximum delivery time $\tau$. If Phase 2 comes after Phase 3, it will become harder to satisfy this constraint. Therefore, in the optimal shipment cycle, we must have Phase 1, 2 and 3 in order. Following that, the optimal shipment policy problem is reduced to determining the length of each service type offered to regular customers (i.e., fast, no and regular service), $t_1$, $t_2$ and $t_3$ respectively. We display the shipment policy structure based on current premium and regular customer demand without knowledge of WoM in Figure \ref{Fig:cycle}. $\lambda_r$ and $\lambda_p$ denote regular and premium demand rate respectively. Note that the implied values for the two quantities of interest are: the length of the shipment cycle, i.e., $T=\sum_{i=1}^3 t_i$, and the maximum inventory level, i.e., $\lambda_pT+\lambda_rt_1$.

\subsubsection{Modeling the Effect of the Membership Fee on Premium Demand}
\label{Sec:Know}

The fixed membership fee, denoted by $F$, is a cost premium customers pay to receive a better service. Premium customers enjoy this privileged service for a period of length $M$. We define $c_1(F)$, a function of membership fee $F$, to represent the potential premium market demand, i.e., maximum premium demand that can be captured by the e-tailer. 
 It is expressed as the product of the number of premium customers denoted by $N(F)$ and demand rate of a premium customer denoted by $\Delta(F)$. When the fee gets higher, becoming a premium member will be less attractive for a given number of purchases per unit time. Therefore, $N(F)$ is a non-increasing function of $F$. On the other hand, each premium customer is expected to order more as the fee increases, i.e., $\Delta(F)$ is non-decreasing, due to the psychological effects of not considering the fixed fee as a sunk cost (see \cite{dick1998impact}, \cite{haller2014sunk}, and \cite{soman2001mental} for discussions in various fields of psychology).

\subsubsection{Modeling the Effect of WoM on Premium Demand}

WoM is a strong means for customers to share information about their service experiences. Customers can share information on various aspects of the service; we call these information \emph{signals}. Signals are reflections of customer perceptions. It has long been known that perception of service quality is determined by comparing the expected and realized service qualities (\cite{parasuraman1990guidelines} and \cite{cardozo1965experimental}). Earlier, we discussed different explanations of how regular customers may perceive the delivery time. The results in the literature suggest that regular customers' perception on the speed of the service is determined by a metric that compares the realized (actual) delivery time and the promised (maximum) delivery time, $\tau$. When the realized service quality is better than the declared service quality, regular customers have a good perception of the service, and provide positive evaluations. This information shared by WoM is available to everyone through online sources. When premium customers receive that information, their perception is simply just the opposite; observing that regular customers are happy with their service experience makes premium customers regret their membership purchase and decreases their incentive to become a member in the future. 

We presume that the signals on the differential between the declared and realized regular service qualities can be measured by a scalar function of the current shipment policy on $[0, 1]$, and denote it by $\theta$, and as $\theta$ decreases, the demand rate for premium membership increases. Different signals can be defined on different aspects of the service. Let $\{\theta_1, \theta_2, 
\dots, \theta_S\}$ be the set of $S$ signals. Observing these signals, a premium customer can develop a perception determined by $\theta=\sum_{i=1}^{S} w_i \theta_i$ where $w_i$ represents the importance of signal $i$, $\sum_{i=1}^{S} w_i=1$ and $w_i\geq 0$ for all $i$.

For the shipment policy presented in Figure \ref{Fig:cycle}, we model two signals that compare the declared and realized service qualities as follows:

\begin{itemize}
\item[] $\theta^{MDT}=\frac{t_3}{\tau}$ as the signal on the \emph{maximum delivery time},
\item[] $\theta^{NPS}=\frac{t_2+t_3}{t_1+t_2+t_3}$ as the signal on the \emph{frequency of providing non-premium service to regular customers},
\end{itemize}

Premium customers' reaction to a received signal $\theta$ is based on their internal response function, which is denoted by $R(\theta)$. We assume non-strategic premium customers, hence premium demand is realized according to the response function on the current signal $\theta$ as follows:
\begin{equation}
\lambda_p \gets R(\theta)
\end{equation}

If customers rely only on the company-released information, then there will be no change in premium demand. On the other hand, if customers are highly sensitive to the information they collect, there can be fluctuations in premium demand. Therefore, the response function should be able to capture the customer sensitivity. We model the response function to an observed signal $\theta \in [0, 1]$ as follows:
\begin{equation}  R(\theta) = c_1(F)\theta^{c_2} \label{Eq:Response}\end{equation}

where $c_1(F)$ is the potential premium market size as defined in Section \ref{Sec:Know}, and $c_2$ is a parameter reflecting premium customer sensitivity to signal $\theta$. As $c_2$ increases, the impact of the signal on premium customer demand increases. $c_2=0$ means completely insensitive premium customers, i.e., WoM has no effect, and $c_2=1$ means linear sensitivity.

In the rest of the section, we study the e-tailer's shipment policy problem in two scenarios: In the first one, the e-tailer does not know how WoM affects demand. On one side, premium customers adjust their service choice observing signals on the service level of regular customers, and on the other side the e-tailer observes the current demand rate and determines the policy accordingly. In the second scenario, we study the case where the e-tailer has knowledge regarding the impact of WoM, and integrates this knowledge into operational decision making.

\subsection{Optimal Shipment Policy without Knowledge on WoM}
\label{Sec:NoKnowWoM}

Without knowledge of WoM, the e-tailer is not aware that signals sent by regular customers affect the perception of premium membership and premium demand as a result. Every time she solves the problem, she treats the most recently observed premium demand as an exogenous parameter.

We assume that the signals on service quality is instantly observable through WoM, and premium customers can react by altering the overall premium demand rate. Similarly, the e-tailer can observe the change in demand and respond with a new shipment policy.

\textbf{Optimal Shipment Policy (given premium and regular customer demand):} Based on current regular and premium demand, we formulate the shipment policy problem to maximize the average profit rate. Our formulation is based on the shipment policy characterization presented earlier in Section \ref{Sec:operations} (Figure \ref{Fig:cycle}). The mathematical model is as follows:
\begin{eqnarray}
  \text{\textbf{\textbf{[M1]}}}\quad  \underset{t_1, t_2, t_3}{\text{max}} & r\lambda_p + \frac{r \lambda_r(t_1+t_3)}{t_1+t_2+t_3}-\frac{h\lambda_p(t_1+t_2+t_3)}{2}-\frac{h\lambda_r t_1^2}{2(t_1+t_2+t_3)}-\frac{K}{t_1+t_2+t_3}
   \label{Eq:ObjNI} \\
      \text{s. t.} & \nonumber\\
       & t_{3} \leq \tau  \label{Eq:Cons1NI} \\
       & 0 \leq t_{1}, t_{2}, t_{3}  \label{Eq:Cons2NI}
\end{eqnarray}

where the objective function in (\ref{Eq:ObjNI}) is the average profit rate, constraint (\ref{Eq:Cons1NI}) ensures that a regular customer waits at most $\tau$ days until she receives her order, and constraints (\ref{Eq:Cons2NI}) ensure non-negativity. In this setting, the company operates on a predetermined membership fee, so it is not a decision variable.

Based on the problem parameters, the objective function can be neither convex or concave, therefore the Karush-Kuhn-Tucker (KKT) conditions are not sufficient for optimality. Therefore, we need to generate the entire set of candidate solutions to solve the problem. However, some solutions can be eliminated. Considering the types of the services offered to regular customers, Phase 3 is the most profitable for the company. So, if the company does not operate in Phase 3, they do not operate in Phases 1 and 2 either. Since we ignore no-operation cases, we ignore the solutions with $t_3=0$. Moreover, based on Lemma \ref{Lem:Obs0}, certain solutions can be identified as dominated by others.

\newtheorem{Lemma}{Lemma}
\begin{Lemma} \label{Lem:Obs0}
When the e-tailer has no knowledge of WoM, the optimal shipment policy has the following properties \vspace{-.1in}
\begin{itemize}
\item[(i)] if $t_1=0$, then $t_2=0$,
\item[(ii)] if $t_3<\tau$, then $t_1=0$.
\end{itemize}
\end{Lemma}

Lemma \ref{Lem:Obs0} follows similarly to Lemma \ref{Lem:Obs} presented in Section \ref{Sec:CompleteInfo}; for brevity, we omit the proof of Lemma \ref{Lem:Obs0}.

Based on Lemma \ref{Lem:Obs0}, there remain only four cases to analyze as shown in Table \ref{Tab:Cases}.

\begin{table}[H]
\scalefont{0.8}
	\begin{center}		
		\caption{Undominated cases for the shipment policy problem with no knowledge on WoM (model \textbf{[M1]}).}		
		\label{Tab:Cases}		
		\begin{tabular}{ c |>{\centering\arraybackslash}p{2.9cm}|>{\centering\arraybackslash}p{2.9cm}|>{\centering\arraybackslash}p{2.9cm} } 			
			\hline \hline			
			 &   Fast service for regular customers& No service for regular customers &  On-time regular delivery\\
			Case &(i.e., $t_{1}> 0$) &(i.e., $t_{2}> 0$) & (i.e., $t_{3}=\tau$) \\ \hline 			
			$i$&$\times$ & $\times$ & $\checkmark$\\ 			
			$ii$ & $\times$  & $\times$  & $\times$ \\			
			$iii$ & $\checkmark$ & $\times$  & $\checkmark$\\			
			$iv$& $\checkmark$& $\checkmark$ & $\checkmark$\\			
			\hline \hline			
		\end{tabular}		
	\end{center}	
	\normalsize	
\end{table}

Case $i$ in Table \ref{Tab:Cases} has a trivial solution, i.e., $t_{1}=t_{2}=0$ and $t_{3}=\tau$. Case $ii$ is the situation where $\tau$ is large enough that it does not limit the company's operations.  In cases $iii$ and $iv$, the e-tailer finds it profitable to send some of the regular orders from the local warehouse even though this incurs additional holding cost. In case $iv$, the company strategically loses some regular demand during the shipment cycle. When we enforce KKT conditions, we find a single solution in each case as shown in Table \ref{Tab:NoInfoSol}. The analysis of the KKT conditions is presented in Appendix \ref{App:NIsolutions}.

\begin{table}[H]
\scalefont{0.8}
\centering
\caption{Undominated solutions for the shipment policy problem with no knowledge on WoM (model \textbf{[M1]}).}
\label{Tab:NIsolutions}
\begin{tabular}{c|c |c| c}
\hline 
\hline
& \multicolumn{3}{c}{Phase Lenghts} \\\cline{2-4}
       \multirow{1}{10mm}{Cases}  &\multirow{1}{15mm}  {\centering $t_{1}$} &   \multirow{1}{15mm}  {\centering $t_{2}$} & \multirow{1}{15mm}  {\centering $t_{3}$}       \\ \hline 
$i$ & 0 & 0 & $\tau$ \\
$ii$ &0 & 0 &$\sqrt{\frac{2K}{h \lambda_{p}}}$ \\
$iii$ & $ \sqrt{\frac{h \lambda_{r} \tau^{2}+2 K}{h (\lambda_{p}+\lambda_{r})}} -\tau$ & 0 & $\tau$	\\
$iv$ & $\frac{r}{h}$ & $\sqrt{\frac{2 hK-2 h \lambda_{r} r \tau-\lambda_{r} r^{2}}{h^2 \lambda_{p}}}-\frac{r}{h}-\tau$  & $\tau$ \\  \hline \hline         
\end{tabular}
\label{Tab:NoInfoSol}
\end{table}

\begin{figure}[H]
\begin{center}
\includegraphics[width=0.65\textwidth]{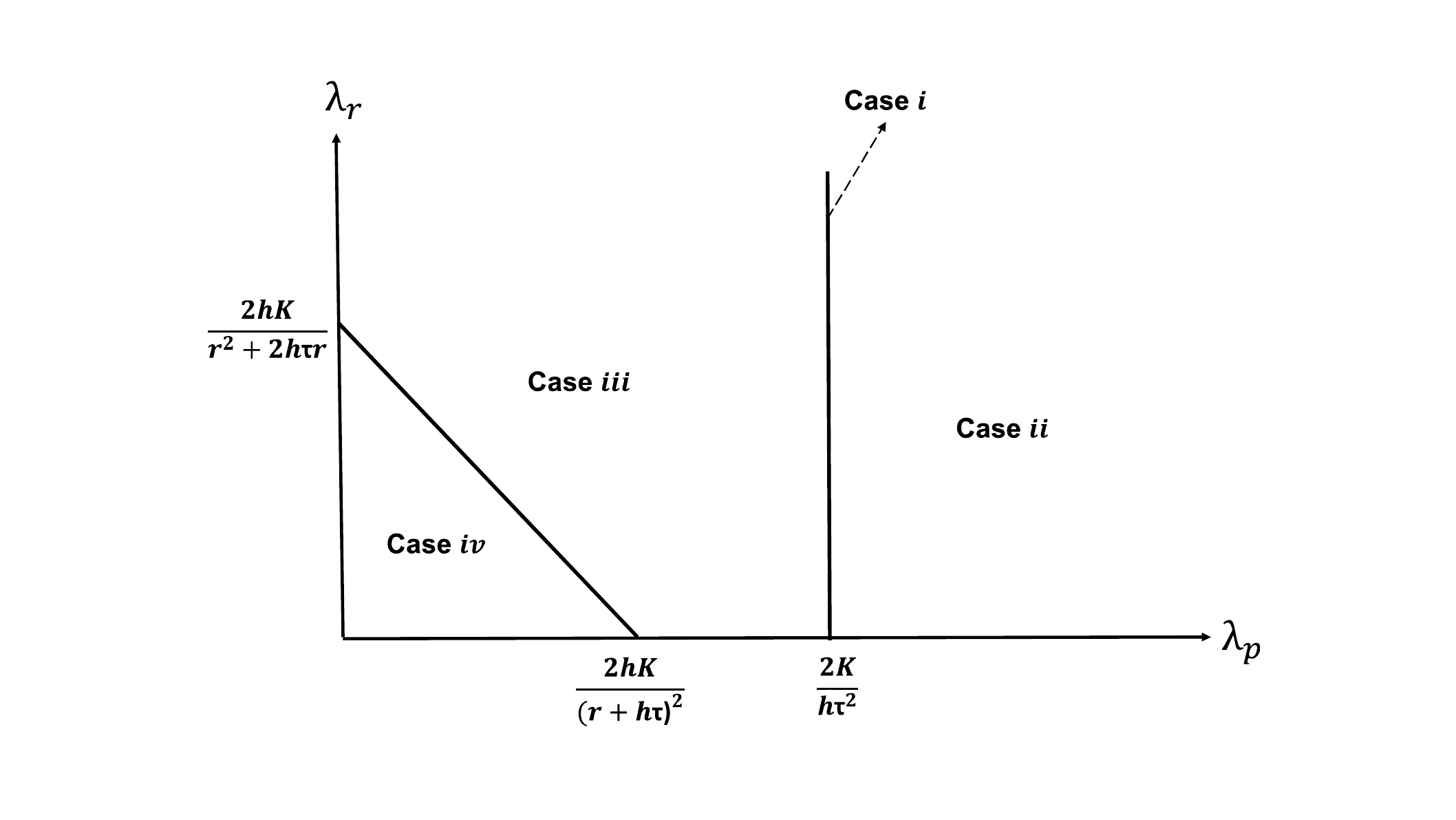}
\caption{Solution regions for the shipment policy problem with no knowledge on WoM (model \textbf{[M1]}).}
\label{Fig:Regions}
\end{center}
\end{figure}

For each solution to be feasible, there are certain conditions that must be satisfied. These conditions determine regions for the demand rates $\lambda_{p}$ and $\lambda_{r}$ expressed as a function of problem parameters. The regions are displayed in Figure \ref{Fig:Regions}.

\textbf{Analysis of the Equilibrium:} Conventionally, the e-tailer would expect the shipment policy determined by solving the problem in model \textbf{[M1]} to be valid in the long run since it optimizes the operations based on the observed demand. However, after customers use the service, they share their experiences creating post-service quality perception which changes customer demand, and decreases the validity of the current shipment policy. As an example, (although this solution is never optimal) if $t_1 > 0$ and $t_2 = t_3 = 0$ in the shipment cycle, there will be no incentive to become a premium member since all regular customers receive their orders immediately. So, there is feedback between the operational decisions and demand. Focusing on this, we analyze the long-term effect of WoM on operations and customer demand.

For the remaining part of this section, we will limit our findings to the maximum delivery time signal, i.e., $\theta=\theta^{MDT}$.

According to the solution space in Figure \ref{Fig:Regions}, if $\frac{2K}{h\tau^2}\geq \lambda_{p}$, then the realized maximum delivery time $t_{3}$ is equal to $\tau$. Since the declared and realized qualities of the regular service are the same, $\theta^{MDT}=1$. However, if $\frac{2K}{h\tau^2}<\lambda_{p}$, the optimal shipment policy contains only Phase 3. In that case, the realized maximum delivery time for regular orders $t_{3}$ is shorter than the declared delivery time $\tau$. Regular customers send positive signals regarding the speed of the regular delivery ($\theta^{MDT}<1$) decreasing the perceived benefit of the premium membership. Based on this post-service information exchange and emerged perceptions, customer choices may shift requiring the company to update the shipment policy according to new demand rates. Theorem \ref{Prop:Eq} characterizes the behavior of premium demand rate when the e-tailer has no knowledge of the impact of WoM.
  
 \newtheorem{Theorem}{Theorem}
\begin{Theorem} \label{Prop:Eq}
	Let $\overline{\lambda_p}$ be the long term premium demand rate under no knowledge on WoM. 
	
	When premium customers respond to the maximum delivery time signal, $\overline{\lambda_p}$ is characterized as follows:
\begin{enumerate}
     \item[] \emph{(i)}  If $c_{1}(F)\leq \frac{2K}{h\tau^2}$, then the company captures the entire potential premium market, i.e., $\overline{\lambda_p}$ converges to $c_{1}(F)$.
    \item[] \emph{(ii)}  If $c_{1}(F)>\frac{2K}{h\tau^2}$ and $c_{2}\geq2$, then $\overline{\lambda_p}$ cycles between  $c_{1}(F)\big{(}\frac{2K}{c_{1}(F)h\tau^{2}}\big{)}^{\frac{c_{2}}{2}}$ and $c_{1}(F)$.
    \item[]  \emph{(iii)} If $c_{1}(F)>\frac{2K}{h\tau^2}$ and $c_{2}<2$, then $\overline{\lambda_p}$ converges linearly to $c_{1}(F)\big{(}\frac{2K}{c_{1}(F)h\tau^{2}}\big{)}^{\frac{c_{2}}{c_{2}+2}}$, and constraint on the maximum delivery time of regular orders becomes redundant thereafter.
\end{enumerate}
\end{Theorem}

The proof is given in Appendix \ref{App:PropProof}. According to Theorem \ref{Prop:Eq}, when the potential premium market size $c_1(F)$ is small, the sensitivity of premium customers does not impact the long run behavior of premium demand. When $c_1(F) \leq \frac{2K}{h\tau^2}$, the operational constraint (\ref{Eq:Cons1NI}) is binding, and the realized service quality is the same as the declared quality for the regular service; hence, the company can capture the entire potential premium market. On the other hand, when $c_1(F)$ is large, regular customers receive a faster service than what it was declared to be, and send positive signals about the regular service. If premium customers are highly sensitive, they show a dramatic reaction to the observed signal resulting in a cyclic demand rate. However, if they are moderately sensitive (i.e., $c_2<2$), then the premium demand rate converges. 

\subsection{Management of Operations and Customer Perception under Knowledge of WoM}
\label{Sec:CompleteInfo}

When the e-tailer is aware of the impacts of WoM, she approaches the shipment policy problem from a different perspective since she can foresee how the current decisions affect the future demand. Under knowledge of WoM, the company focuses on determining a shipment policy that nudges premium demand into what is the most profitable in the long run. The long term behavior of premium demand also depends on the potential premium market size ($c_1(F)$) determined by $F$; hence, membership fee also emerges as a natural decision variable.

We formulate the shipment policy problem under knowledge of WoM in model \textbf{[M2]}. The model is built under two assumptions: (1) the company uses the  shipment cycle structure characterized in Section \ref{Sec:operations}, and (2) the company aims to achieve stationary customer demand through WoM for ease of operations.

\begin{eqnarray}
\text{\textbf{[M2]}} \quad \underset{t_1, t_3, T, \lambda_p,F, \theta}{\text{max }} &r\lambda_p+\frac{r\lambda_r(t_1+t_3)}{T}
 +\frac{F \lambda_p}{\Delta(F)M}-\frac{K}{T}-\frac{hT\lambda_p}{2}-\frac{h\lambda_r t_1^2}{2T} &\label{Eq:ObjCI} \\
\text{s. t.} & \nonumber\\
& t_{3} \leq \tau  \label{Eq:Cons1CI} \\
& t_1+t_{3} \leq T  \label{Eq:Cons2CI} \\
 &\lambda_p= R(\theta)\label{Eq:Cons3CI} \\
& f_{min} \leq F\leq f_{max} \label{Eq:Cons4CI} \\
& 0 \leq t_1, t_3  \label{Eq:Cons5CI} 
\end{eqnarray}

where $t_1$, $t_3$, $T$ and $F$ are the primary decision variables, and $\theta$ is the observed WoM signal. Notice that $\lambda_p$ is an auxiliary decision variable; we do not determine it directly, but it is the equilibrium premium demand rate realized based on the membership fee and shipment policy choice. Objective function in (\ref{Eq:ObjCI}) is the average profit rate. Constraints (\ref{Eq:Cons1CI}) and (\ref{Eq:Cons2CI}) follow from the previous model. Constraint (\ref{Eq:Cons3CI}) computes the equilibrium premium demand rate using the assumption of stationary demand. Constraints (\ref{Eq:Cons4CI}) and (\ref{Eq:Cons5CI}) provide bounds for the decision variables.

Management of customer perceptions is one of the two main focuses of model \textbf{[M2]}. The perception of regular service quality depends on the realized and declared service qualities. While the former is the output of the model (decision of the company), the latter is imposed by market dynamics. As an example, consider a case where $\theta$ is a function of the maximum delivery time signal only, i.e., $\theta=\frac{t_3}{\tau}$. We can interpret the parameter $\tau$ as a reflection of market competition; $\tau$ is smaller in competitive markets. If the market is highly competitive, the company cannot operate with long lead times. Hence, $\tau$ is considered to be an exogenous parameter. Whenever market dynamics change, the company can update the shipment policy so that the signal is managed as intended. For example, the e-tailer can respond to changes in $\tau$ by changing $t_3$ so that the perceptions are still realized in her favor. However, this cannot always be done easily due to the other aspects of the decision making, e.g., corresponding costs and operational constraints. 

In the rest of this section, focusing on the tradeoff between the two drivers of the problem, i.e., managing customer perceptions and maximizing profit, we partially characterize of the optimal shipment policy under knowledge of WoM.

\textbf{Properties of the Optimal Solution for \textbf{[M2]}}: In our analysis, we assume that $N(F)$  is a linear decreasing function of $F$, $\Delta(F)$ is a constant function, and premium customers are linearly sensitive, i.e.,  $N(F)=a-bF$, $\Delta(F)=\delta$, and $c_{2}=1$. We also assume that $\theta=\theta^{MDT}$. We consider other scenarios in our computational experiments.

We first show that a similar observation in Lemma \ref{Lem:Obs0} also holds for the optimal shipment policy problem under knowledge of WoM. 

\begin{Lemma} \label{Lem:Obs}
When the e-tailer has knowledge of WoM, and premium customers respond to the maximum delivery time signal, the company does not operate in Phase 2 if they are not operating in Phase 1, and does not operate in Phase 1 when $t_3<\tau$.
\end{Lemma}

We prove Lemma \ref{Lem:Obs} in Appendix \ref{App:LemmaObs}. 
Based on Lemma \ref{Lem:Obs}, if the operational constraint in (\ref{Eq:Cons1CI}) is not binding, then $t_1=t_2=0$, $t_3<\tau$ and $T=t_3$ in the optimal solution. In that case, a greater $\tau$ relaxes a constraint that is already redundant. So, we would expect the optimal shipment policy to remain the same. However, Proposition \ref{Prop:tau_t3} below shows that this is not the case. Moreover, Proposition \ref{Prop:tau_z} shows that profit rate decreases under the same condition.

 \newtheorem{Proposition}{Proposition}
\begin{Proposition}
	\label{Prop:tau_t3}
Consider that the e-tailer has knowledge of WoM and premium customers respond to the maximum delivery time signal.
		
When $t_3<\tau$ in the optimal solution, greater $\tau$ results in offering a worse service to regular customers, i.e., $t_3$ increases.
\end{Proposition}

We prove Proposition \ref{Prop:tau_t3} in Appendix  \ref{App:tau_t3proof}. If one presumes $t_3$ to remain the same, greater $\tau$ values result regular customers to expect a longer wait and hence develop a better perception of the service. As a result, regular customers send more positive signals, which would hurt the premium demand. To prevent this, the company increases the realized maximum delivery time. This is an interesting result in terms of showing how a company (that can actually deliver the orders earlier) is forced to deliver the orders later to be able to manage customers' perceptions and so the premium demand. Therefore, while we would expect an increase in profit by relaxing the operational constraints, we show that this is not always the case due to the associated difficulties of managing customer perceptions in the proposition below. 
 \begin{Proposition}
 	\label{Prop:tau_z}
 	
Consider that the e-tailer has knowledge of WoM and premium customers respond to the maximum delivery time signal.

When $t_3<\tau$ in the optimal shipment cycle, the profit decreases as $\tau$ gets greater.
 \end{Proposition}
 
The proof is in Appendix \ref{App:tau_zproof}. The impact of an increase in $\tau$ on the profit is a natural result of how $t_3$ changes with $\tau$. If $t_3$ is currently less than $\tau$ in the shipment cycle, then the company is not limited by the operational constraints. A greater $\tau$ is expected not to worsen the solution. However, it forces the company to operate with a longer Phase 3 due to Proposition \ref{Prop:tau_t3}. Hence, it becomes more costly to manage customer perceptions when the declared service quality is worse. 

These results reveal the counter-intuitive impacts of the changes in problem parameters when customers' choices are affected by their perceptions. Such cases occur because the two drivers of the problem may work in different directions: while the profit maximization considers a known and static demand rate, management of customer perception focuses on creating a more profitable demand through signaling which can be done by offering a seemingly suboptimal solution for a static demand rate. When the second force is stronger than the first one, the solution differs from our prevalent expectations.

The effects of WoM are especially significant when the operational constraints are not binding. When $\tau$ is small (which can be interpreted as a sign of a competitive market), the feasible region is small; the company's decisions are already limited. There is not much room for an action to manage perceptions. On the other hand, when $\tau$ is greater, not only the feasible region is larger but also managing the perceptions becomes more challenging. Therefore, problem parameters affect the extend of the impact of WoM. 

\textbf{What Happens If We Have Partial Information of WoM?} While the complete knowledge of WoM is required to make optimal shipment policy decisions considering customer perceptions, there are cases where the partial information on WoM would still imitate the full knowledge. Such cases occur when the operational constraints are likely to be binding.  Proposition \ref{Prop:recoverability} characterizes one of those cases.

 \begin{Proposition}
 	\label{Prop:recoverability}
 	Consider that the e-tailer has knowledge of WoM, and premium customers are linearly sensitive to the maximum delivery time signal, i.e., $c_2=1$.

 	Let $\overline{\lambda_p}$ be the long-term premium demand rate without knowledge of WoM, and $\lambda_p^{Eq}$ be the equilibrium premium demand rate achieved under knowledge of WoM.

When an e-tailer does not have knowledge of WoM and is exogenously provided with only the optimal fee $F$ determined by model \textbf{[M2]},  $\overline{\lambda_p}=\lambda_p^{Eq}$ only if $c_1(F) \leq \frac{2K}{h\tau^2}$. Otherwise, $\overline{\lambda_p} \leq \lambda_p^{Eq}$. \end{Proposition}

The proof is given in Appendix \ref{App:Recover}. The proposition implies that in competitive markets with small $\tau$, partial information of the optimal shipment policy determined under knowledge of WoM can mimic knowledge of WoM. In such cases, once the optimal fee $F$ is computed, both shipment policy models (\textbf{[M1]} and \textbf{[M2]}), can produce the same results. On the other hand, for greater $\tau$'s, this cannot happen, and the shipment policy must be determined by considering the impact of WoM on customer demand. 

Our theoretical results show counter-intuitive impacts of WoM on operations and demand in different aspects. In the following section, we present numerical demonstrations for these results and also reveal new observations through computational experiments.

\section{Computational Studies}
\label{Sec:Comp}

In our computational studies, we performed experiments using different WoM signals, i.e., MDT and NPS, and different types of relationships between the membership fee and the number of potential premium customers, namely $N(F)=a-bF$ and $N(F)=alog(b-F)$. We performed these experiments on a rich set of parameters as shown in Table \ref{Tab:Param}; however, for brevity we present only a small subset to accompany our discussion in the paper.

\begin{table}[H]
\scalefont{0.8}
	\begin{center}		
		\caption{Parameters for the computational experiments.}		
		\label{Tab:Param}	
\begin{tabular}{r c l}
		\hline
 $r$ &$\in$&   $\{8, 16, 48\}$\\
 $K$ & $\in$ &      $\{2000, 3000, 4000\}$\\
 $\tau$ & $\in$ &      $\{1, 1.5, 2, 3, 4, 5, 6, 7\}$\\
 $c_{2}$   & $\in$& $\{0.1, 0.2, 0.5, 1, 2, 3\}$\\
  $\delta$ & $\in$ &      $\{0.56, 5\}$\\
 $\lambda_r$ &=&$50$\\
 $h$ &=& 4\\
 $(a, b)$ &=&$(100, 1),$ \text{ for linear $c_1(F)$}\\
$(a, b)$ &=&$(20, 101),$ \text{ for logarithmic $c_1(F)$}\\ 
  $M$ & $\in$ &      $\{Monthly, Lifetime\}$\\
\hline 
\end{tabular}
\end{center}
\end{table}

Our experiments are composed of two parts. In the first part, we determine the optimal solutions for different scenarios under knowledge of WoM using model \textbf{[M2]}, and discuss the impact of WoM on the optimal shipment policy and demand. In the second part, we input the optimal fee with the corresponding equilibrium premium demand rate determined by \textbf{[M2]} to the shipment policy problem with no knowledge on WoM, which is modeled in \textbf{[M1]}, and analyze the behavior of the shipment policy and demand in the long run. The purpose of the second experiment is to implement our Proposition \ref{Prop:recoverability} numerically to see when the partial information of WoM can imitate knowledge of WoM. Using these computational results, we also test the two assumptions we used while building model \textbf{[M2]}.

We display the main numerical results for the first experimental setting in Tables \ref{Tab:MDTlinearc1}, \ref{Tab:MDTlogc1}, \ref{Tab:NPSlinearc1} and \ref{Tab:NPSlogc1} respectively for MDT signal with linear and logarithmic $N(F)$ functions, and NPS signal with linear and logarithmic $N(F)$ functions. In those tables, columns $t_1$, $t_2$, $t_3$ and $F$ show the values of the optimal solutions for the decision variables. $\lambda_p$ shows the corresponding equilibrium premium demand rate, and $\pi$ is the optimal profit rate. No-WoM Decision corresponds to the second experiment. When we input the optimal solutions of the first experiment to model \textbf{[M1]}, and simulate the long term behavior, three cases may occur: we either maintain the initial equilibrium solution (denoted by Opt-Eq), deviate from the initial equilibrium and converge to a non-optimal equilibrium (denoted by Non-opt-Eq), or observe a cyclic shipment policy and customer demand (denoted by Cycles).


\begin{table}[H]
\scalefont{0.85}
\centering
\caption{Optimal shipment policy under knowledge of WoM - MDT signal, $N(F)=a -bF$, $r=8$, $h=4$, $K=2000$, $\lambda_r=50$, $a=100$,  $b=1$, $\delta=5$, $M=30$, $f_{min}=10$,  $f_{max}=100$.}
\label{Tab:MDTlinearc1}
\begin{tabular}{c |c |c c c c c c| c}
\hline \hline 
$c_2$ &{$\tau$}     & $t_1$ & $t_2$ & $t_3$ & $F$   & $\lambda_{p}$ & $\pi$ &No-WoM Decision \\
\hline
\multirow{5}{*}{1}&1.0	&	0.45	&	0.00	&	1.00	&	10	&	450.00	&	1331.73	&	Opt-Eq	\\	
&1.5	&		0.00	&	0.00	&	1.50	&	10	&	450.00	&	1346.67	&	Non-opt-Eq	\\	
&2.0	&		0.00	&	0.00	&	2.00	&	10	&	450.00	&	1230.00	&	Non-opt-Eq	\\	
&5.0	   &	0.00	&	0.00	&	2.75	&	10	&	247.58	&	307.98	&	Non-opt-Eq	\\	
&6.0	&		0.00	&	0.00	&	2.84	&	10	&	213.15	&	204.14	&	Non-opt-Eq	
\\ \hline
\multirow{5}{*}{3} & 1.0		&	0.45	&	0.00	&	1.00	&	10	&	450.00	&	1331.73	&	Opt-Eq	\\	
&1.5		&	0.00	&	0.00	&	1.50	&	10	&	450.00	&	1346.67	&	Cycles	\\	
&2.0		&	0.00	&	0.00	&	2.00	&	10	&	450.00	&	1230.00	&	Cycles	\\	
&5.0	&	1.71	&	0.00	&	5.00	&	100	&	0.00	&	58.36	&	Opt-Eq	\\	
&6.0		&	1.48	&	0.00	&	6.00	&	100	&	0.00	&	103.34	&	Opt-Eq	\\	
\hline \hline 
\end{tabular}
\end{table}

\begin{table}[H]
\scalefont{0.85}
\centering
\caption{Optimal shipment policy under knowledge of WoM - MDT signal, $N(F)=a log(b-F)$, $r=8$, $h=4$, $K=2000$, $\lambda_r=50$, $a=20$,  $b=101$, $\delta=5$, $M=30$, $f_{min}=10$,  $f_{max}=100$.}
\label{Tab:MDTlogc1}
\begin{tabular}{c |c |c c c c c c| c}
\hline \hline 
{$c_2$} & {$\tau$}      & $t_1$ & $t_2$ & $t_3$ & $F$   & $\lambda_{p}$ & $\pi $ &No-WoM Decision \\
\hline
\multirow{6}{*}{1} & 1.0		&	0.05	&	0.00	&	1.00	&	10.00	&	451.09	&	2138.87	&	Opt-Eq	\\	
&1.5	&		0.00	&	0.00	&	1.50	&	10.00	&	451.09	&	2018.84	&	Non-opt-Eq	\\	
&2.0	&		0.00	&	0.00	&	2.00	&	10.00	&	451.09	&	1734.42	&	Non-opt-Eq	\\	
&5.0	&	0.00	&	0.00	&	2.47	&	10.00	&	222.89	&	691.88	&	Non-opt-Eq	\\	
& 6.0	&	0.00	&	0.00	&	2.53	&	10.00	&	190.54	&	576.64	&	Non-opt-Eq	\\	
\hline
\multirow{6}{*}{3}&1.0		&	0.05	&	0.00	&	1.00	&	10.00	&	451.09	&	2138.87	&	Opt-Eq	\\	
&1.5	&		0.00	&	0.00	&	1.50	&	10.00	&	451.09	&	2018.84	&	Cycles	\\	
&2.0	&		0.00	&	0.00	&	2.00	&	10.00	&	451.09	&	1734.42	&	Cycles	\\	
&5.0		&	0.00	&	0.00	&	3.40	&	45.21	&	126.77	&	295.74	&	Cycles *	\\	
&6.0	&		0.78	&	0.00	&	6.00	&	100.00	&	0.00	&	243.53	&	Opt-Eq	\\	
\hline \hline 
\end{tabular}
\end{table}



\begin{table}[H]
\scalefont{0.85}
\centering
\caption{Optimal shipment policy under knowledge of WoM - NPS signal, $N(F)=a -bF$, $\tau=1$, $h=4$, $\lambda_r=50$, $a=100$,  $b=1$, $\delta=5$, $M=30$, $f_{min}=10$,  $f_{max}=100$.}
\label{Tab:NPSlinearc1}
\begin{tabular}{c| c| c  |c c c c c c |c}
\hline \hline 
$K$ & $r$ &{$c_2$} & $t_1$ & $t_2$ & $t_3$ & $F$   & $\lambda_{p}$ & $\pi$ & No-WoM Decision\\ \hline 
\multirow{6}{*}{3000}	&	\multirow{3}{*}{16}	&	0.10	&	0.35	&	0.23	&	1.00	&	10.00	&	438.76	&	4441.47	&	Non-opt-Eq	\\	
	&		&	0.20	&	0.00	&	0.57	&	1.00	&	10.00	&	450.00	&	4415.80	&	Non-opt-Eq	\\	
	&		&	1.00	&	0.00	&	0.56	&	1.00	&	10.00	&	450.00	&	4415.75	&	Non-opt-Eq	\\	\cline{2-10}
	&	\multirow{3}{*}{48}	&	0.10	&	0.01	&	0.00	&	1.00	&	10.00	&	449.69	&	20130.16	&	Non-opt-Eq	\\	
	&		&	0.20	&	0.00	&	0.00	&	1.00	&	10.00	&	450.00	&	20130.00	&	Non-opt-Eq	\\	
	&		&	1.00	&	0.00	&	0.00	&	1.00	&	10.00	&	450.00	&	20130.06	&	Non-opt-Eq	\\
\hline	\hline
\multirow{6}{*}{4000}	&	\multirow{3}{*}{16}			&	0.10	&	0.45	&	0.45	&	1.00	&	10.00	&	437.94	&	3867.46	&	Non-opt-Eq	\\	
	&		&	0.20	&	0.00	&	0.89	&	1.00	&	10.00	&	450.00	&	3835.91	&	Non-opt-Eq	\\	
	&		&	1.00	&	0.00	&	0.89	&	1.00	&	10.00	&	450.00	&	3835.89	&	Non-opt-Eq	\\\cline{2-10}
	&	\multirow{3}{*}{48}	&	0.10	&	0.20	&	0.15	&	1.00	&	10.00	&	442.87	&	19257.85	&	Non-opt-Eq	\\	
	&		&	0.20	&	0.00	&	0.33	&	1.00	&	10.00	&	450.00	&	19230.10	&	Non-opt-Eq	\\	
	&		&	1.00	&	0.00	&	0.33	&	1.00	&	10.00	&	450.00	&	19230.08	&	Non-opt-Eq	\\	
\hline \hline
\end{tabular}
\end{table}


\begin{table}[H]
\scalefont{0.85}
\centering
\caption{Optimal shipment policy under knowledge of WoM - NPS signal, $N(F)=a log(b-F)$,$\tau=1$, $h=4$, $\lambda_r=50$, $a=100$,  $b=1$, $\delta=5$, $M=30$, $f_{min}=10$,  $f_{max}=100$.}
\label{Tab:NPSlogc1}
\begin{tabular}{c| c| c  |c c c c c c| c}
\hline \hline 
$K$ & $r$ &{$c_2$} & $t_1$ & $t_2$ & $t_3$ & $F$   & $\lambda_{p}$ & $\pi$  &No-WoM Decision\\ \hline 
\multirow{6}{*}{3000}	&	\multirow{3}{*}{16}&	0.10	&	0.34	&	0.23	&	1.00	&	10.00	&	440.16	&	4455.13	&	Non-opt-Eq			\\	
	&		&	0.20	&	0.00	&	0.56	&	1.00	&	10.00	&	451.09	&	4429.85	&	Non-opt-Eq			\\	
	&		&	1.00	&	0.00	&	0.56	&	1.00	&	10.00	&	451.09	&	4429.84	&	Non-opt-Eq			\\\cline{2-10}
	&		\multirow{3}{*}{48}	&	0.10	&	0.01	&	0.00	&	1.00	&	10.00	&	450.52	&	20180.12	&	Non-opt-Eq			\\	
	&		&	0.20	&	0.00	&	0.00	&	1.00	&	10.00	&	451.09	&	20180.14	&	Non-opt-Eq			\\	
	&		&	1.00	&	0.00	&	0.00	&	1.00	&	10.00	&	451.09	&	20180.19	&	Non-opt-Eq			\\	
\hline\hline
\multirow{6}{*}{4000}	&	\multirow{3}{*}{16}	&	0.10	&	0.46	&	0.45	&	1.00	&	10.00	&	438.84	&	3880.43	&	Non-opt-Eq			\\	
&		&	0.20	&	0.00	&	0.88	&	1.00	&	10.00	&	451.09	&	3849.26	&	Non-opt-Eq			\\	
&		&	1.00	&	0.00	&	0.88	&	1.00	&	10.00	&	451.09	&	3849.24	&	Non-opt-Eq			\\\cline{2-10}	
	&	\multirow{3}{*}{48}	&	0.10	&	0.19	&	0.15	&	1.00	&	10.00	&	444.14	&	19306.37	&	Non-opt-Eq			\\	
	&		&	0.20	&	0.00	&	0.33	&	1.00	&	10.00	&	451.09	&	19279.51	&	Non-opt-Eq			\\	
	&		&	1.00	&	0.00	&	0.33	&	1.00	&	10.00	&	451.09	&	19279.46	&	Non-opt-Eq			\\	
\hline \hline 
\end{tabular}
\end{table}


We discuss our findings from the computational experiments below.

\subsection{Impact of exogenous parameters}

In our problem setting, exogenous parameters are maximum delivery duration imposed by market dynamics $\tau$, premium customer sensitivity $c_2$, membership duration $M$, holding cost $h$, revenue per item $r$, and shipment cost $K$. Among those $\tau$, $M$ and $c_2$ directly determine the impact of WoM on customer demand. The impact of changes in $h$, $r$ and $K$ are what we would expect in a classical shipment policy problem, e.g., as $K$ increases shipment cycle length increases, as $h$ increases the e-tailer is less willing to serve regular customers in Phase 1. Since we are mainly interested in the impact of WoM, we limit our discussion to $\tau$, $M$ and $c_2$.

\textbf{Impact of $\tau$: }In Proposition \ref{Prop:tau_t3}, we show that under knowledge of WoM, when $t_3<\tau$ in the optimal solution, a greater $\tau$ increases $t_3$. Our computations numerically demonstrate this on different scenarios. For example, in Tables \ref{Tab:MDTlinearc1} and \ref{Tab:MDTlogc1}, we observe that when $\tau=5$, the constraint $t_3\leq \tau$ is not binding for $c_2=1$. However, when $\tau$ becomes 6, $t_3$ also increases. The reason is that as $\tau$ increases, regular customers expect a slower delivery. If the company still operates according to the previous delivery speed, regular customers become relatively happier which may hurt premium customer demand. To avoid this, the e-tailer is forced to offer a worse service to regular customers when $\tau$ is greater.

Our computational results also illustrate that increasing $\tau$, i.e., relaxing an operational constraint, may hurt the profit as shown in Proposition \ref{Prop:tau_z}. For example, in Table \ref{Tab:MDTlinearc1}, when $\tau=1$, the operational constraint is binding, i.e., $t_3=\tau$ for all $c_2$ values. In that case, increasing $\tau$ to 1.5 improves the profit. However, as $\tau$ gets larger, the profit decreases because managing customer perception becomes financially more burdensome under maximum delivery time signal. Therefore, although an operational relaxation is not expected to cause a decrease in the profit, it can do so under knowledge of WoM.

Responding to the changes in $\tau$ becomes more difficult as it gets greater. In Table \ref{Tab:MDTlinearc1}, we also observe the impact of large $\tau$'s on operations under the maximum delivery time signal. When $\tau$ exceeds a certain threshold, managing customer perception profitably becomes so difficult that the e-tailer prefers to remove the premium service as a response. To achieve this, the  e-tailer increases the membership fee to its maximum value so that no customer finds it attractive to become a premium member. We can see such a case when $\tau=6$ and $c_2=3$. This threshold decreases when premium customers are more sensitive; while the threshold is 7 when $c_2=1$, it is 6 when $c_2=2$, and 5 when $c_2=3$. Once the threshold is exceeded, increasing $\tau$ improves the profits as we would expect in a system where WoM has no impact. Moreover, once $\tau$ exceeds that threshold, the e-tailer increases the shipment cycle length by serving some regular customers in Phase 1, i.e., providing premium service to regular customers. This is an example of a case where a system which is designed to offer different types of services becomes very challenging to manage due to WoM, and the company has to remove some service types.

\textbf{Impact of lifetime membership: }The decision of removing premium customers from the system depends on other factors as well. One of these is the length of the membership term. As the duration of the membership increases, the additional benefit of a premium customer decreases. In such cases, dropping the premium service from the system is easier especially when premium and regular demand rates are comparable. In Table \ref{Tab:MDTlinearlifetime}, we observe such a case. Notice that in this data set, the potential premium demand and regular demand rates are equal. When the duration of membership is 30 days, the e-tailer keeps the premium service until $\tau$ increases to 7, on the other hand for large membership durations, the e-tailer removes the premium service even when $\tau$ is 4.

\begin{table}[H]
\scalefont{0.85}
\centering
\caption{Optimal shipment policy under knowledge of WoM - MDT signal, $N(F)=a-bF$, $r=8$, $h=4$, $c_2=1$, $K=1000$, $\lambda_r=50$, $a=20$,  $b=101$, $f_{min}=10$,  $f_{max}=100$, $\delta=0.56$.}
\label{Tab:MDTlinearlifetime}
\begin{tabular}{c| c  |c c c c c c}
\hline \hline 
$M$ &{$\tau$}     & $t_1$ & $t_2$ & $t_3$ & $F$   & $\lambda_{p}$ & $\pi $ \\
\hline
	\multirow{5}{*}{Monthly Membership}	&	1.0	&	1.44	&	0.00	&	1.00	&	24.04	&	42.20	&	97.72		\\	
	&	1.5	&	1.10	&	0.00	&	1.50	&	26.63	&	40.76	&	148.11		\\	
	&	2.0	&	0.81	&	0.00	&	2.00	&	30.12	&	38.82	&	183.34		\\	
&	5.0	&	0.08	&	0.00	&	5.00	&	68.03	&	17.76	&	237.16	\\	
&	6.0	&	0.47	&	0.00	&	6.00	&	91.10	&	4.95	&	244.63	\\	
\hline
	\multirow{5}{*}{Lifetime Membership}	&	1.0	&	1.35	&	0.00	&	1.00	&	10.00	&	50.00	&	61.92		\\	
	&	1.5	&	0.98	&	0.00	&	1.50	&	10.00	&	50.00	&	110.05		\\	
	&	2.0	&	0.65	&	0.00	&	2.00	&	10.00	&	50.00	&	141.70		\\	
	&	5.0	&	0.92	&	0.00	&	5.00	&	100.00	&	0.00	&	216.78	\\	
	&	6.0	&	0.78	&	0.00	&	6.00	&	100.00	&	0.00	&	243.53	\\	
\hline \hline
\end{tabular}
\end{table}

\textbf{Impact of customer sensitivity: }When premium customers are more sensitive, the reactions to the observed signal are bigger. As a result, the e-tailer is forced to take stronger measures to manage customer perceptions. For example, in Tables \ref{Tab:MDTlinearc1} and \ref{Tab:MDTlogc1}, we observe a greater increase in $t_3$ and greater decrease in profit as $\tau$ gets larger under MDT signal when premium customers are more sensitive.


\subsection{Impact of partial information on WoM}

\textbf{Recovering the optimal solution to model [M2] through the partial information on WoM: }In the second set of experiments, we observe whether the e-tailer that is not aware of the WoM communication can recover the optimal solution to model \textbf{[M2]} when exogenously provided with the optimal membership fee and the corresponding equilibrium premium demand. As indicated in the last columns of Tables \ref{Tab:NPSlinearc1} and \ref{Tab:NPSlogc1}, under NPS signal, the e-tailer cannot recover the optimal solution through the partial information since NPS signal changes the structure of the optimal solution. Although this is not the case for MDT signal, recovering knowledge of WoM using the partial information is still rare.

\begin{figure}[H]
	\begin{center}
		\includegraphics[width=0.6\textwidth]{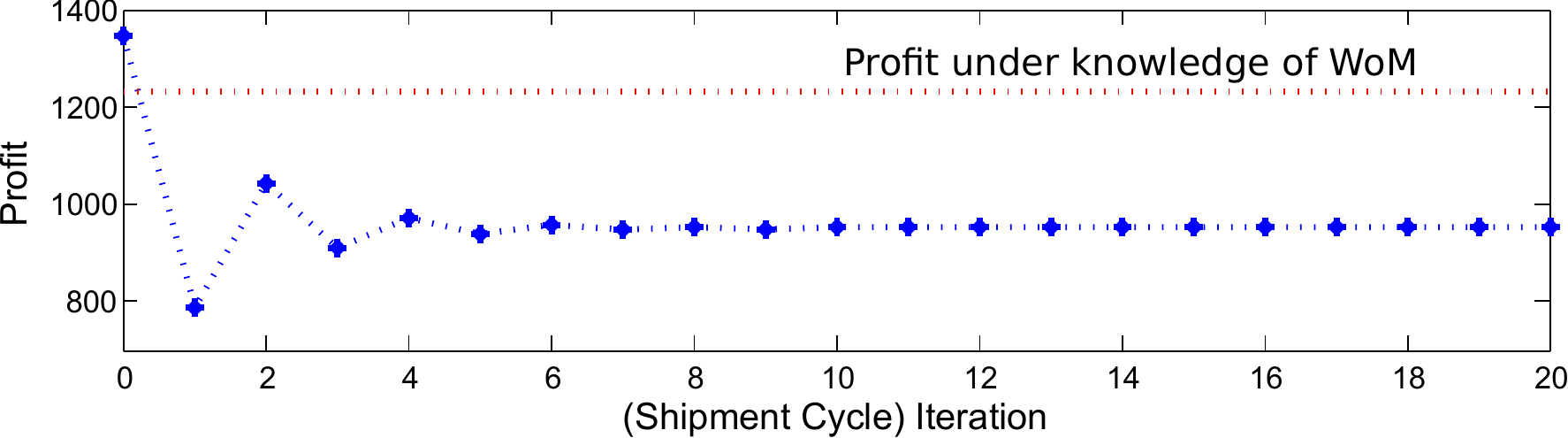}
		\caption{Profit without knowledge of WoM - MDT signal, $N(F)=a -bF$, $r=8$, $h=4$, $K=2000$, $\lambda_{R}=50$, $a=100$,  $b=1$, $\delta=5$, $M=30$, $\tau=2$, $c_2=1$, $F=10$.}
		\label{Fig:MDTNon-opt-Eq}
	\end{center}
\end{figure}

For demonstration of the second experiment, we present the results for two parameter sets under MDT signal from Table \ref{Tab:MDTlinearc1}. These are $\tau=2$ and $c_2=1$, and $\tau=2$ and $c_2=3$. For these instances, we let model \textbf{[M1]} start from the equilibrium premium demand rate determined under knowledge of WoM, and compute the shipment policy iteratively using the newly realized premium demand. Table \ref{Tab:MDTNon-opt-Eq} shows the corresponding shipment policy decisions and the resulting premium demand rate for $\tau=2$ and $c_2=1$. Figure \ref{Fig:MDTNon-opt-Eq} displays the profit rates where the horizontal dashed line shows the corresponding optimal profit rate under knowledge of WoM.

\begin{table}[H]
\centering
\small
\caption{Shipment policy and corresponding premium demand without knowledge of WoM - MDT signal, $N(F)=a -bF$, $r=8$, $h=4$, $K=2000$, $\lambda_{R}=50$, $a=100$,  $b=1$, $\delta=5$, $M=30$, $\tau=2$, $c_2=1$, $F=10$.}
\label{Tab:MDTNon-opt-Eq}
\begin{tabular}{c|c c c c c c c c c c c }\hline \hline
$Iteration$	&	0	&	1	&	2	&	3	&	4	&	5	&	6	&	7	&	8	&	9	&	10\\	\hline
$\lambda_{p}$	&	450.00	&	335.41	&	388.50	&	360.98	&	374.49	&	367.67	&	371.07	&	369.37	&	370.22	&	369.79	&	370.00\\	
$t_{1}$	&	0.00	&	0.00	&0.00	&0.00	&	0.00	&	0.00	&	0.00	&	0.00	&	0.00	&	0.00	&	0.00\\	
$t_{2}$	&	0.00	&	0.00	&	0.00	&	0.00	&	0.00	&	0.00	&	0.00	&	0.00	&	0.00	&	0.00	&	0.00\\	
$t_{3}$	&	1.49	&	1.73	&	1.60	&	1.66	&	1.63	&	1.65	&	1.64	&	1.65	&	1.64	&	1.64	&	1.64\\	
\hline \hline 
\end{tabular}
\end{table}

In Figure \ref{Fig:MDTNon-opt-Eq}, we observe that using model \textbf{[M1]} the e-tailer can achieve a higher profit than the average profit that can be achieved under knowledge of WoM in the first iteration. The reason is that she ignores the impact of the current decisions on customer demand, and myopically focuses on profits. However, the resulting shipment policy leads to a sharp decrease in premium demand in the next period, and the e-tailer's response is to worsen the regular service quality as seen in Table \ref{Tab:MDTNon-opt-Eq}. The e-tailer's policy is only reactive to the changes in demand without knowing that those changes are endogenous to the operational decisions. In the long run, premium demand converges to a non-optimal equilibrium, and the company earns 29\% less compared to integrating knowledge of WoM in operational decision making.

\begin{figure}[H]
	\begin{center}
		\includegraphics[width=0.6\textwidth]{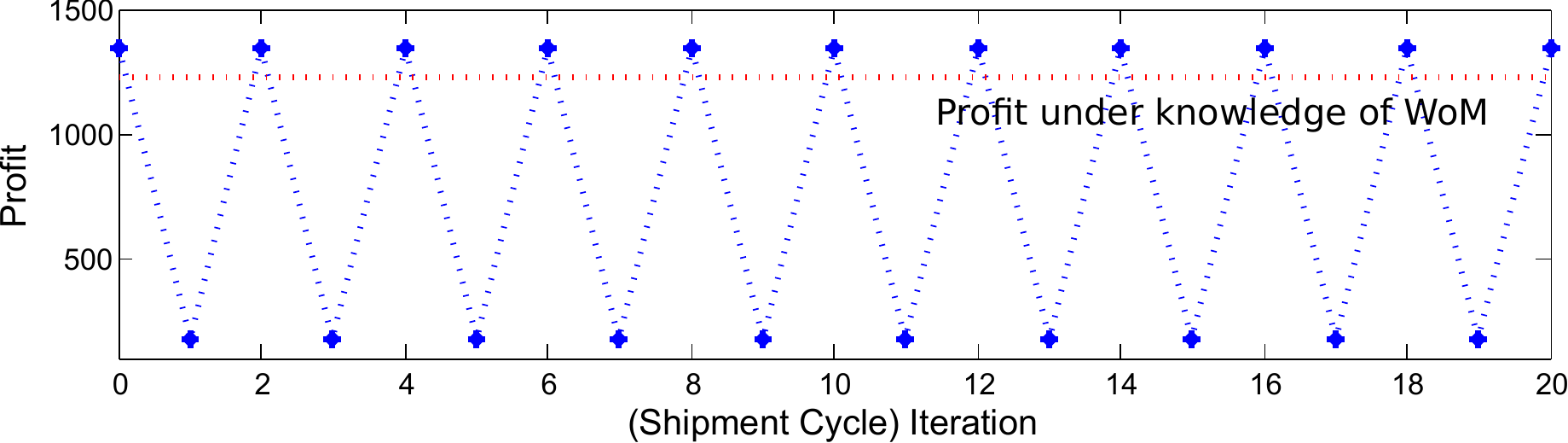}
		\caption{Profit without knowledge of WoM - MDT signal, $N(F)=a -bF$, $r=8$, $h=4$, $K=2000$, $\lambda_{R}=50$, $a=100$,  $b=1$, $\delta=5$, $M=30$, $\tau=2$, $c_2=3$, $F=10$.}
		\label{Fig:MDTcycles}
	\end{center}
\end{figure}

We present the same set of results for $\tau=2$ and $c_2=3$ in Table \ref{Tab:MDTcycles} and Figure \ref{Fig:MDTcycles}. We make similar observations to the previous case: ignoring the impact WoM and focusing only on the current profit, higher profits can be achieved in the first iteration; yet, in the long run the average profit is smaller. However, this time we observe a cyclic premium demand and shipment policy mostly due to higher sensitivity of premium customers. When $c_2=3$, the reaction of premium customers to the observed signal becomes much bigger, and it prevents convergence. 

\begin{table}[H]
\centering
\small
\caption{Shipment policy and corresponding premium demand without knowledge of WoM - MDT signal, $N(F)=a -bF$, $r=8$, $h=4$, $K=2000$, $\lambda_{R}=50$, $a=100$,  $b=1$, $\delta=5$, $M=30$, $\tau=2$, $c_2=3$, $F=10$.}
\label{Tab:MDTcycles}
\begin{tabular}{c|c c c c c c c c c c c }\hline \hline
$Iteration$	&	0	&	1	&	2	&	3	&	4	&	5	&	6	&	7	&	8	&	9	&	10\\	\hline
$\lambda_{p}$	&	450.00	&	186.34	&	450.00	&	186.34	&	450.00	&	186.34	&	450.00	&	186.34	&	450.00	&	186.34	&	450.00\\	
$t_{1}$	&	0.00	&	0.25	&	0.00	&	0.25	&	0.00	&	0.25	&	0.00	&	0.25	&	0.00	&	0.25	&	0.00\\	
$t_{2}$	&	0.00	&	0.00	&	0.00	&	0.00	&	0.00	&	0.00	&	0.00	&	0.00	&	0.00	&	0.00	&	0.00\\	
$t_{3}$	&	1.49	&	2.00	&	1.49	&	2.00	&	1.49	&	2.00	&	1.49	&	2.00	&	1.49	&	2.00	&	1.49\\	
\hline \hline 
\end{tabular}
\end{table}

\textbf{Cyclic shipment policies: }Another interesting observation from the computational results corresponds to parameters $\tau=5$ and $c_2=3$ marked by (*) in Table \ref{Tab:MDTlogc1}. In this case, when we simulate the shipment policy decisions without knowledge of WoM according to the second experiment, the shipment policy follows a cyclic behavior, and the average profit is higher than the profit attained under knowledge of WoM in the first experiment. The reason is that model \textbf{[M2]} limits the solution to stationary policies when we integrate knowledge of WoM into our shipment problem (assumption 2 for model \textbf{[M2]}). The assumption is due to the fact that stationary policies provide easier management, and provide a consistent perception of the service quality rather than a cyclic reputation. Although it is very rare that cyclic policies perform better than stationary policies in our experiments (only two cases are observed), these results suggest that cyclic policies may also be considered when it is appropriate.

\subsection{Impact of WoM signal type}
We earlier discussed that when solving the shipment policy problem without knowledge of WoM, the e-tailer does not serve in Phase 2 if not serving in Phase 1 due to the associated profitability of those phases. This property is still sustained with MDT signal. Therefore, under knowledge of WoM and MDT signal, all solutions obey Lemma \ref{Lem:Obs}. However, this is not the case for NPS signal (i.e., when premium customers care about whether regular customer orders are delivered fast or not). Recall the two drivers of the objective function in model \textbf{[M2]}: profit maximization for a static demand rate and management of customer perception to attain a more profitable demand rate. While the former obeys Lemma \ref{Lem:Obs}, the latter may act in a different manner. When premium customers are sensitive to NPS signal, perception management driver in the problem tries to prevent serving in Phase 1. Therefore, we observe cases where $t_2>0$ while $t_1=0$ in the optimal solution if the second force is stronger than the first one.

Looking at the results in Tables \ref{Tab:NPSlinearc1} and \ref{Tab:NPSlogc1} on NPS signal, we observe that the solutions obey Lemma \ref{Lem:Obs} when premium customer sensitivity is very low, e.g., $c_2=0.1$. In that case, the reactions of premium customers to the observed NPS signal is very small that it does not hurt premium demand much to provide fast service to some regular customers. However, as $c_2$ increases, this is not the case anymore. In the optimal solution, fast service to regular customers becomes less frequent and it finally disappears (measured by $t_1$), while the e-tailer increases the frequency of losing regular customers (measured by $t_2$). Moreover, we also observe that deviation of the optimal shipment policy from the behavior characterized in Lemma \ref{Lem:Obs} is larger for larger $c_2$ values.

\section{Conclusion}
\label{Sec:Conc}

Motivated by the behavioral literature that discusses the impact of WoM on customer perceptions, in this paper we study the interaction between operational decisions and demand in an e-tailer system. We consider an e-tailer that offers two types of services, namely premium and regular, and customers whose choices depend on the perceived differential between these service types. Building on the WoM literature, we model a demand structure endogenous to operational decisions for premium customers who may regret if the perceived difference in service quality is not as large as expected. 

We show that the long term behavior of customer choices depends on factors such as premium customers' sensitivity, potential premium market size and the information provided on the maximum delivery time. When the e-tailer has no knowledge of WoM, she has no control on customer perceptions, and simply reacts to changes in customer demand. In the long run, customer demand may converge to an undesired point, and the company is stuck at a suboptimal operational policy. Alternatively, in markets with a large potential for premium demand and sensitive premium customers, demand may follow a cyclic behavior resulting in cyclic operational policies. This creates a belief that the company's service quality is not consistent. None of these results can be predicted ignoring the impact of WoM on perceptions. Considering how information exchange between customers affect their perception and demand for a service, we analyze the impact of WoM in an operational decision making.  We also conduct an extensive computational study to gain further insight. The managerial insights are explained in detail in the previous section. Here we would like to emphasize the ones that challenge the conventional wisdom and support them through industry examples.

First, we show that beating the market competition or extensive marketing efforts can have adverse effects under operational restrictions. Such actions may result in an undesired (or excessive) customer switching due to positive WoM. If the operational capacity is not ready to respond to the increase in demand, the company cannot maintain the high service quality perception, and the company reputation is hurt. This insight partially explains some bankruptcies occurred at the beginning of 2000's. When EToys founded in 1997, they wanted to compete with companies such as Amazon, Toys R Us and Wal-Mart, and followed an aggressive marketing policy (\cite{EToys}). As a result, they increased their customer base quickly. However, they failed to deliver orders on time as their operational capacity was not ready to handle  large demand. The decreased reputation and their efforts to invest in infrastructure to respond to high customer demand were among the main reasons preparing their bankruptcy in 2001. A similar phenomenon was also observed with Charming Charlie bankruptcy in 2017. Having founded in 2004, they showed a rapid growth offering a wide range of options, and faced a widespread consumer demand. This made it difficult for them to understand the dynamics affecting the demand and respond to them. The company reported that reducing the scale of their operations provided them the opportunity to better read and react to what customers want (\cite{CharmingCharlie}). From these examples, we see that outperforming the market standards can result negatively for companies especially if they are not operationally ready to respond to large demand inflow. Moreover, it makes it difficult to understand dynamics (possibly endogenous to company operations and reputation) affecting customer demand.

We also show that relaxations in operational constraints may not always benefit the company; on the contrary, they may be hurtful if these constraints are associated with customer perceptions. Under the knowledge of WoM, there are two drivers affecting operational decisions, one is profit maximization and the other one is managing customer perceptions. While a relaxation in an operational constraint enlarges the solution space and so is desirable from the first driver's perspective, it can also interfere with how customer perceptions are realized and introduce additional financial burden to manage those perceptions. As an example, in the paper we show that in less competitive markets, the e-tailer can declare a longer delivery duration for regular customers, which relaxes an operational constraint; however, this can decrease the profit due to the impact of WoM.

Finally, seeking a stationary policy (although it is highly favored in the literature due to its ease of management) may result in suboptimal solutions. Our results suggest that companies may also consider changing their service quality from time to time to maintain a desired customer perception and demand for their services. Such strategy is of course subject to whether the underlying application is suitable for cyclic policies or not.

The insights of this paper can be applied to other types of businesses where different services/products are offered, and customers make their choices based on the perceived differences. Our model manages customer perceptions under the assumption that the e-tailer has no control over the competition in the market. Future research can focus on leveraging the knowledge of WoM for an industry-leading company that can set standards in the market.

\bibliographystyle{chicago}
\bibliography{eWoM}

\appendix
\section{Appendices}

\subsection{Solutions without Knowledge of WoM}
\label{App:NIsolutions}

The Lagrangian relaxation of the  equivalent minimization problem of the shipment policy problem in model \textbf{[M1]} is as follows:

\begin{eqnarray}
L(t_{1}, t_{2}, t_{3}, y)&=&\frac{1}{2} h \lambda_{p}(t_{1}+t_{2}+t_{3})+\frac{h \lambda_{r} t_{1}^{2}}{2 (t_{1}+t_{2}+t_{3})}+\frac{K}{t_{1}+t_{2}+t_{3}}+\frac{\lambda_{r} r t_{2}}{t_{1}+t_{2}+t_{3}} \nonumber \\ &&-y_1 (\tau-t_{3})- y_{2}t_{1}-y_{3}t_{2}-y_{4}t_{3},
\end{eqnarray}

where $y=(y_1, y_2, y_3, y_4)$ is the vector of corresponding nonnegative dual variables for the constraints in (\ref{Eq:Cons1NI}) and (\ref{Eq:Cons2NI}).


The solutions to the problem must satisfy the following KKT conditions.

\begin{itemize}
\item Stationary: $ \frac{\partial L(t_{1}, t_{2}, t_{3},y)} {\partial t_{i}}=0,   \quad \forall i $
\item Complementary slackness: $y_{1}(\tau-t_{3})=0$, $y_{2}t_{1}=0$, $y_{3}t_{2}=0$, $y_{4}t_{3}=0$
\item Primal feasibility: $0\leq t_{3} \leq \tau$, $0\leq t_{1}$, $0\leq t_{2}$
\item Dual feasibility: $y_{i} \geq 0$, $\forall i$
\end{itemize}

We solve for the solutions that satisfy the KKT conditions for the undominated cases shown in Table \ref{Tab:Cases}, and present the results in Table \ref{Tab:NIsolutions}. 
\subsection{Proof of Theorem \ref{Prop:Eq}}
\label{App:PropProof}
 \begin{proof}

When premium customers respond to maximum delivery time signal, future premium demand is realized according to equation (\ref{Eq:Response}).

Based on the solution to the shipment policy problem in model \textbf{[M1]}, post-service maximum delivery time signal $\theta^{MDT}$ can be expressed as a function of the current premium demand rate $\lambda_{p}$ as follows:
	
	\begin{eqnarray}
	\label{Eq:mdtcases}
	\theta^{MDT}=\begin{cases}
	1,      \quad\quad \textrm{if } \lambda_{p} \leq \frac{2K}{h\tau^2},\\
	\frac{\sqrt{\frac{2K}{h\lambda_{p}}}}{\tau}, \quad \textrm{otherwise.}
	\end{cases}
	\end{eqnarray}

Let us denote the premium demand rate prior to solving the shipment policy problem for the $k^{th}$ time as $\lambda_{p}^{k}$.

First, consider that $c_1(F)\leq \frac{2K}{h\tau^2}$. In that case $\lambda_p^{k} \leq \frac{2K}{h\tau^2}$ for all $k$ values. According to equation (\ref{Eq:mdtcases}), all future MDT signals will be realized as 1, and premium demand rate immediately converges to $c_1(F)$.

Now, consider that $c_1(F) >\frac{2K}{h\tau^2}$. Let $\lambda_p^{0}$ be the current premium demand rate. 

If $\lambda_p^{0}\leq \frac{2K}{h\tau^2}$, then $\lambda_p^1=c_1(F)$. Following this, $\lambda_p^{2}$ is realized as

\begin{equation}
\lambda_p^{2}=c_1(F) \bigg{(}\frac{2K}{h\tau^2c_1(F)}\bigg{)}^{\frac{c_{2}}{2}}.
\end{equation}

For ease of notation, let us denote $\frac{2K}{h\tau^2c_1(F)}$ by $w$, and note that $w<1$.

If $c_2\geq 2$, then $\lambda_p^{2}\leq\frac{2K}{h\tau^2}$, and $\lambda_p^3=c_1(F)$. Therefore, premium demand rate cycles between $c_1(F)$ and $c_1(F) \bigg{(}\frac{2K}{h\tau^2c_1(F)}\bigg{)}^{\frac{c_{2}}{2}}$.

If $c_2<2$, then $\lambda_p^{2}>\frac{2K}{h\tau^2}$, and $\lambda_p^3$ is realized as 

\begin{equation}
\lambda_p^{3}=c_1(F) \bigg{(}\frac{2K}{h\tau^2c_1(F) w^{\frac{c_2}{2}}}\bigg{)}^{\frac{c_{2}}{2}}.
\end{equation}

To prove the convergence of premium demand rate when $c_1(F)>\frac{2K}{h\tau^2}$ and $c_2<2$, first we will show that $\lambda_p^{k}$ can be written as $c_1(F)w^{B_k}$  where $B_k<1$ for all $k\geq 2$ by induction. As the base case, we have $\lambda_p^2=c_1(F)w^{\frac{c_2}{2}}$ where $c_2<2$. Assume $\lambda_p^{k}=c_1(F)w^{B_k}$ where $B_k<1$ for some $k>2$. Since $w<1$, $\lambda_p^{k}>c_1(F)w = \frac{2K}{h\tau^2}$. By equation \ref{Eq:mdtcases}, $\lambda_p^{k+1}$ is realized as 

\begin{equation}
\lambda_p^{k+1}=c_1(F)  \bigg{(} \frac{2K}{h\tau^2 c_1(F) w^{B_k}}\bigg{)}^{\frac{c_2}{2}}.
\end{equation}

So, $\lambda_p^{k+1}=c_1(F)w^{B_{k+1}}$, where $B_{k+1}=(1-B_k)\frac{c_2}{2}<1$. By induction, $\lambda_p^{k}=c_1(F)w^{B_k}$ where $B_k<1$ for all $k\geq2$. 

Since $w<1$, $\frac{2K}{h\tau^2}=c_1(F)w < c_1(F)w^{B_{k}}$ for all $k\geq2$. Hence, once the premium demand exceeds $\frac{2K}{h\tau^2}$, it remains so when $c_1(F)>\frac{2K}{h\tau^2}$ and $c_2<2$. Therefore, $\lambda_p^{0}>\frac{2K}{h\tau^2}$ case does not need to be analyzed. 

Based on the induction, $\lambda_p^{k}$ can be written as follows.
\begin{equation}
\label{Eq:lambdapk}
\lambda_{p}^{k} = c_1(F)w^{\sum_{i=1}^{k-1}\big{(}\frac{c_2}{2}\big{)}^{i}(-1)^{i-1}}.
\end{equation}

As $k\rightarrow \infty$, $\sum_{i=1}^{k-1}\big{(}\frac{c_2}{2}\big{)}^{i}(-1)^{i-1}$ can be expressed as an alternating series that converges to $\frac{c_2}{c_2+2}$. Therefore, premium demand rate $\lambda_p$ converges to $c_{1}(F)\big{(}\frac{2K}{c_{1}(F)h\tau^{2}}\big{)}^{\frac{c_{2}}{c_{2}+2}}$ when $c_2<2$ and $c_1(F)> \frac{2K}{h\tau^2}$.
%
%
%
%
%
\end{proof}
\subsection{Proof of Lemma \ref{Lem:Obs}}
\label{App:LemmaObs} 
 
 \begin{proof}
We prove Lemma \ref{Lem:Obs} by showing that a solution to the shipment policy problem in model \textbf{[M2]} under the knowledge on maximum delivery time signal cannot be improved by the following operations:
 
 \begin{itemize}
 	\item[(i)] decreasing the length of Phase 1 by $\epsilon$ and increasing the length of Phase 2 by $\epsilon$,
 	\item[(ii)] decreasing the length of Phase 3 by $\epsilon$ and increasing the length of Phase 1 by $\epsilon$,
 \end{itemize}

where $\epsilon>0$.

Let $(t_1, t_3, T, F)$ be a feasible solution. To prove (i), we show that the solution $(t_1-\epsilon, t_3, T, F)$ has a worse objective function value than $(t_1, t_3, T, F)$. Let $z$ and $z_{\epsilon}$ be the objective function values corresponding to $(t_1, t_3, T, F)$ and $(t_1-\epsilon, t_3, T, F)$ respectively. $z-z_{\epsilon}$ is equal to $\epsilon(r-ht_1+h\epsilon/2$). Therefore, $z$ is always greater than $z_\epsilon$ if $t_1\leq \frac{r}{h}$. Now, we show that this condition always holds.

The Lagrangian function of the equivalent minimization problem of the shipment policy problem in model \textbf{[M2]} is as follows:

\begin{eqnarray}
\label{Eq:Lagfull}
L(t_{1}, t_{3}, T, F, \omega)&= &\frac{K}{T}+\frac{hT}{2}(a-bF)\frac{\delta t_3}{\tau}+\frac{h\lambda_r t_1^2}{2T}-r(a-bF)\frac{\delta t_3}{\tau}- \frac{r \lambda_r(t_1+t_3)}{T}  \nonumber \\
&&-\frac{F(a-bF)t_3}{\tau M}-\omega_1t_1-\omega_2t_3-\omega_3(\tau-t_3)-\omega_4(T-t_1-t_3) \nonumber \\
&&-\omega_5(F-f_{min})-\omega_6(f_{max}-F), \label{Eq:M2Lag}
\end{eqnarray}

where $\omega=(\omega_1, \omega_2, ... , \omega_6)$ is the vector of corresponding nonnegative dual variables for the constraints in (\ref{Eq:Cons1CI})-(\ref{Eq:Cons5CI}). The following is one of the equalities that must be satisfied for stationarity:

\begin{equation}
\label{Eq:fullKKT_t1}
\frac{\partial L(t_{1}, t_{3}, T, F, \omega)}{ \partial t_1}=\frac{h \lambda_r t_1}{T}-\frac{r \lambda_r }{T}-\omega_1 +\omega_4=0
\end{equation}

We are interested in finding a bound on $t_1$, basically we are interested in the cases where $t_1>0$. In that case, the corresponding dual variable $\omega_1$ is 0. Based on equality (\ref{Eq:fullKKT_t1}), we have that $t_1\leq \frac{r}{h}$ since $\omega_4 \geq 0$. This completes the first part of the proof.

To prove (ii), we follow a similar procedure and show that the solution $(t_1+\epsilon, t_3-\epsilon, T, F)$ has a worse objective function value than $(t_1, t_3, T, F)$. Now, let $z$ and $z_{\epsilon}$ be the objective function values corresponding to $(t_1, t_3, T, F)$ and $(t_1+\epsilon, t_3-\epsilon, T, F)$ respectively. $z-z_{\epsilon}$ is computed as follows:

\begin{equation}
\label{Eq:z-zepsilon2}
\frac{(a-bF)\epsilon}{\tau} (\delta r -\frac{\delta h T}{2}+\frac{F}{M}) +\frac{h \lambda_r}{2T}(2t_1\epsilon +\epsilon^2).
\end{equation}

We need to show that the expression in (\ref{Eq:z-zepsilon2}) is always nonnegative. We again use KKT conditions for that. One of the stationarity conditions that must be satisfied is as follows:

\begin{equation}
\frac{{\partial L(t_1, t_3, T, F, \omega) } }{\partial t_3} =\frac{hT}{2}(a-bF)\frac{\delta}{\tau}-\frac{r(a-bF)\delta}{\tau}-\frac{r \lambda_r}{T}-\frac{F(a-bF)}{\tau M}-\omega_2+\omega_3+\omega_4=0
\label{Eq:fullKKT_t3}
 \end{equation}

By subtracting the expression in (\ref{Eq:fullKKT_t1}) from (\ref{Eq:fullKKT_t3}), we find that

\begin{equation}
\label{Eq:KKT_t3-t1}
\frac{(a-bF)}{\tau}(\delta r-\frac{h\delta}{2}+\frac{F}{M})+\frac{h \lambda_r t_1}{T}=\omega_2+\omega_3.
\end{equation}

The right hand side of equation (\ref{Eq:KKT_t3-t1}) is nonnegative. It follows that (\ref{Eq:z-zepsilon2}) is positive when $\epsilon >0$, which completes the second part of the proof.
 
\end{proof}


 \subsection{Proof of Proposition \ref{Prop:tau_t3}}
 \label{App:tau_t3proof}
 \begin{proof}
 	The proof of Proposition \ref{Prop:tau_t3} is simple but tedious; we explicitly determine the optimal $t_3$ in model \textbf{[M2]} and show that it is increasing in $\tau$. The proposition considers the cases where $t_3< \tau$ in the optimal solution. In this case, by Lemma \ref{Lem:Obs}, we know that $t_1=t_2=0$ in the optimal solution. When we analyze the KKT conditions for this case, we find that the optimal solution is determined by one of the following cases:
 	
 	\begin{itemize}
 	\item[(i)] $f_{min} <F < f_{max}$
 	\begin{equation}
 	t_3=\frac{2(a+b\delta M r)}{3b\delta M r}+\frac{4\cdot2^{1/3}(a+b\delta Mr)^2}{3b\delta M h X}+\frac{X}{3\cdot2^{1/3}b\delta M r},
 	\end{equation}

where $X=\bigg{(}Y+108b^3\delta ^2 M^3 h^2 K\tau+\sqrt{-256(a+b\delta M r)^6+(Y+108b^3\delta^2 M^3 h^2 K \tau)^2}\bigg{)}^{1/3}$, and $Y=16(a+b\delta Mr)^3$.

The partial derivative of $t_3$ with respect to $\tau$ is $\frac{\partial t_3}{\partial \tau}=\frac{-4\cdot 2^{1/3}(a+b\delta M r)^2}{3b\delta Mh X^2}\frac{\partial X}{\partial \tau}+\frac{1}{3\cdot2^{1/3}b\delta M h}\frac{\partial X}{\partial \tau}$, which is always nonnegative since $X\geq 2\cdot 2^{1/3}(a+b\delta Mr)$ and $\frac{\partial X}{\partial \tau}>0$.
 	\item[(ii)] $F=f_{min}$ or $f_{max}$
 	\end{itemize}
 	
 	\begin{equation}
 	t_3=\frac{1}{6 \delta M h (a-b F)} (2 (a-b F) (\delta  M r+F)+ \frac{2\cdot  2 ^{1/3} (a-b F)^2 (\delta  M r+F)^2}{(X+Y)^{1/3}}) + 2^{2/3}(2X+Y)^{1/3}
 	\end{equation}
 	
 	where 
 	$$X=(a-b F)^2 \left(2 (a-b F) (\delta  M r+F)^3+27 \delta ^2 \text{DD}^3 h^2 K\tau\right), $$
 	$$Y=3 \sqrt{3} \sqrt{\delta ^2 M^3 h^2 K \tau (a-b F)^4 \left(4 a (\delta  M r+F)^3-4 b F (\delta  M r+F)^3+27 \delta ^2 M^3 h^2 K\tau\right)}.$$
 	
 	The partial derivative of $t_3$ with respect to $\tau$ has only one negative term, i.e., $-2(a-bF)^2(F+\beta M r)^2$. It also has a positive term of the same magnitude. Therefore, $\frac{\partial X}{\partial \tau}>0$ in this case as well.
 	\end{proof}
 	
  \subsection{Proof of Proposition \ref{Prop:tau_z}}
 \label{App:tau_zproof}
 \begin{proof}
 	Let $(t_1, t_3, T, F)$ be the optimal solution for model \textbf{[M2]} when $\theta=\theta^{MDT}$. Due to Lemma \ref{Lem:Obs}, we know that if $t_3 <\tau$, then $t_1=t_2=0$ and $T=t_3$. So, the optimal objective function value $z$ can be written as follows
 	
 	\begin{equation}
 	\frac{hT(a-bF)\delta t_3}{2\tau^2}+\frac{-r(a-bF)\delta t_3}{\tau^2}+\frac{h \lambda_r t_1^2}{2\tau^2}- \frac{F(a-bF)t_3}{\tau^2M},
 	\end{equation}

and the first derivative of $z$ with respect to $\tau$ is equal to

\begin{equation}
\label{Eq:zpartt3}
\frac{\partial z}{\partial \tau}=\frac{(a-bF)t_3}{\tau^2}\bigg{(} \frac{hT\delta}{2}-r\delta-\frac{F}{M} \bigg{)}.
\end{equation}

Recall that equation (\ref{Eq:M2Lag}) provides the Lagrangian function of the equivalent minimization problem for model \textbf{[M2]}. We know that in the optimal solution $\frac{\partial L(t_{1}, t_{3}, T, F, \omega)  }{\partial t_1}+\frac{\partial L(t_{1}, t_{3}, T, F,\omega)  }{\partial t_3}$ is equal to 0. Since we are considering the case where $0<t_3<\tau$, $\omega_2$ and $\omega_3$ are both 0. Therefore, the summation can be written as follows:

\begin{equation}
\frac{a-bF}{\tau}\bigg{(}\frac{hT\delta}{2}-r\delta -\frac{F}{M}\bigg{)} +\omega_1=0,
\end{equation}

which implies that $T\leq \frac{2r}{h}+\frac{2F}{h\delta M}$ since $\omega_1\geq 0$ and ${a-bF}\geq 0$. Therefore, the derivative of the optimal objective function value $z$ with respect to $\tau$ in Equation (\ref{Eq:zpartt3}) is negative. 
\end{proof}	

\subsection{Proof of Proposition \ref{Prop:recoverability}}

\label{App:Recover}
\begin{proof}
 When customers are linearly sensitive to the observed demand  ($c_2=1$), only cases (i) and (iii) of Theorem \ref{Prop:Eq} are relevant. 
 
To prove the first part of the proposition, we show that if $c_1(F)\leq \frac{2K}{h\tau^2}$, then $t_3=\tau$ in the optimal solution under knowledge of WoM. 

The optimal solution to model \textbf{[M2]} must satisfy equation (\ref{Eq:fullKKT_t1}) and the following stationarity condition.
 
 \begin{equation}
 \label{Eq:Tpartial}
 \frac{\partial L(t_{1}, t_{3}, T, F, \omega) }{ \partial T}=-\frac{K}{T^2}+\frac{h c_1(F)t_3}{2\tau}-\frac{h\lambda_r t_1^{2}}{2T^2}+\frac{r \lambda_r t_1}{T^2}+\frac{r \lambda_r t_3}{T^2}-\omega_4=0.
 \end{equation}
 
Consider the case where $c_1(F)\leq \frac{2K}{h\tau^2}$ in the optimal solution. Assume that the optimal $t_3<\tau$. In that case, $t_1$ and $t_2$ must be 0, and $T=t_3$. So, equation (\ref{Eq:Tpartial}) becomes

\begin{equation}
\frac{-K}{T^2} +\frac{h c_1(F) T}{2\tau}+\frac{r\lambda_r}{T}-\omega_4=0.
\end{equation}

Summing equations (\ref{Eq:fullKKT_t1}) and (\ref{Eq:Tpartial}), we obtain $\frac{-K}{T^2}+\frac{hc_1(F)T}{2\tau}=\omega_1$. Since $T<\tau$, we have $c_1(F)<\frac{2K}{hT^2}$ implying that $\omega_1<0$, which cannot occur. Therefore, when $c_1(F)\leq \frac{2K}{h\tau^2}$ in the optimal solution, then $t_3=\tau$; hence, and both $\overline{\lambda_p}$ and $\lambda_p^{Eq}$ are equal to $c_1(F)$.
 
Now, we prove the second part of the proposition. Since we already showed that $c_1(F)\leq\frac{2K}{h\tau^2}$ implies $t_3=\tau$, $t_3<\tau$ implies that $c_1(F)> \frac{2K}{h\tau^2}$ in the optimal solution.

When $t_3<\tau$, $t_1=t_2=0$  and $T=t_3$ by Lemma \ref{Lem:Obs}. In that case, equations (\ref{Eq:fullKKT_t1}) and (\ref{Eq:Tpartial}) can be written as follows respectively.

\begin{equation}
\label{Eq:newt3}
\frac{r\lambda_r}{T}-\omega_1+\omega_4=0
\end{equation}

\begin{equation}
\label{Eq:newT}
\frac{-K}{T^2}+\frac{hc_1(F)T}{2\tau}+\frac{r \lambda_r}{T}+\omega_4=0
\end{equation}

Summing (\ref{Eq:newt3}) and (\ref{Eq:newT}), we obtain

\begin{equation}
\frac{-K}{T^2}+\frac{hc_1(F)T}{2\tau}=\omega_1,
\end{equation}

where the left hand side is nonnegative since  $\omega_1\geq 0$. This implies $\frac{T^3}{\tau}\geq \frac{2K}{hc_1}$, and $\frac{t_3}{\tau}\geq (\frac{2K}{hc_1(F)\tau})^{1/3}$. Therefore, when $c_1(F)>\frac{2K}{h\tau^2}$, $\lambda_p^{Eq}=c_1(F)\frac{t_3}{\tau}$ is greater than or equal to $\overline{\lambda_p}=c_1(F)(\frac{2K}{hc_1(F)\tau})^{1/3}$.

\end{proof}

\end{document}